\newif\ifspringer
\springerfalse

\ifspringer

\RequirePackage{fix-cm}
\RequirePackage[colorlinks,citecolor=blue,urlcolor=blue]{hyperref}

\documentclass[smallextended]{svjour3}       
\smartqed  
\journalname{4OR}

\else
\documentclass[10pt,a4paper]{article}

\fi

\usepackage[normalem]{ulem}
\usepackage{soul}
\usepackage{amsmath}
\usepackage{amsfonts}
\usepackage{amssymb}         
\usepackage{amsopn}
\usepackage{theorem}          
\usepackage{latexsym}
\usepackage{graphicx}        
\usepackage{mathrsfs}		
\usepackage{psfrag}
\usepackage{algorithmic}
\usepackage{algorithm}
\usepackage{xcolor}
\usepackage{verbatim}
\usepackage{url}
\usepackage{enumerate}
\usepackage{wrapfig}
\usepackage{comment}
\usepackage{tikz}
\usepackage{circuitikz}
\usetikzlibrary{shapes.geometric}
\usetikzlibrary{arrows.meta,arrows}
\usepackage[title]{appendix}

\definecolor{midgrey}{gray}{0.6}
\definecolor{darkgrey}{gray}{0.4}

\ifspringer
\else
\fi

\theoremstyle{changebreak}                
\ifspringer
\newtheorem{thm}[result]{Theorem}
\newtheorem{defn}[result]{Definition}
\newtheorem{lem}[result]{Lemma}
\newtheorem{cor}[result]{Corollary}

\newtheorem{eg}[result]{Example}
\else

 {{\nopagebreak\hspace*{\fill}$\Box$\par\vspace{12pt}}}

\fi

\newcommand{\transpose}[1]{{#1}^{\top}}
\newcommand{\hermitian}[1]{{#1}^{\mathsf{H}}}

\definecolor{plgreen}{rgb}{0,0.5,0}

\newcommand{\re}[1]{{#1}^{\mbox{\sf\scriptsize r}}}
\newcommand{\im}[1]{{#1}^{\mbox{\sf\scriptsize c}}}

\newcommand{\conj}[1]{\,{#1}\!{}^\ast}
\newcommand{\bluem}[1]{{\color{blue}{#1}}}
\newcommand{\redem}[1]{{\color{red}{#1}}}
\renewcommand{\vec}[1]{\overrightarrow{#1}}
\newcommand{\cev}[1]{\overleftarrow{#1}}

\newcommand{\trace}[1]{{\mathsf{Tr}}{#1}}
\newcommand{\seplin}{\noindent\rule{0.5\textwidth}{0.3pt}\linebreak\par}

\newcommand{\LEO}[1]{{\color{black}#1}}

\begin{document}

\ifspringer
\else

\thispagestyle{empty}
\begin{center} 

{\LARGE Mathematical Programming formulations for the Alternating Current Optimal Power Flow problem}
\par \bigskip
{\sc Dan Bienstock${}^1$, Mauro Escobar${}^2$, Claudio Gentile${}^3$, Leo Liberti${}^2$}
\par \bigskip
 {\small
   \begin{enumerate}
   \item {\it IEOR, Columbia University, New York, USA} \\ Email:\url{dano@columbia.edu}
   \item {\it LIX CNRS, \'Ecole Polytechnique, Institut Polytechnique de Paris, F-91128 Palaiseau, France} \\ Email:\url{mauro.escobar-santoro@polytechnique.edu,liberti@lix.polytechnique.fr}
   \item {\it IASI, CNR, Rome, Italy} \\ Email:\url{claudio.gentile@iasi.cnr.it}
   \end{enumerate}
 }
\par \medskip \today
\end{center}
\par \bigskip

\fi

\begin{abstract}
\LEO{Power flow refers to the injection of power on the lines of an electrical grid, so that all the injections at the nodes form a consistent flow within the network. Optimality, in this setting, is usually intended as the minimization of the cost of generating power. Current can either be direct or alternating: while the former yields approximate linear programming formulations, the latter yields formulations of a much more interesting sort: namely, nonconvex nonlinear programs in complex numbers. In this technical survey, we derive formulation variants and relaxations of the alternating current optimal power flow problem.}
  \ifspringer
\else
  \\ {\bf Keywords:} ACOPF, Smart grid, complex numbers.
  \fi
\end{abstract}

\tableofcontents

\section{Introduction}
\label{s:intro}
This paper is about Mathematical Programming (MP) formulations for the Alternating Current Optimal Power Flow (ACOPF) problem. The ACOPF an optimization problem aiming at generating and transporting electrical power at minimum cost.

The basic entity in electricity is the {\it charge}: by analogy, charge is to the electromagnetic force what mass is to gravity. Charge is transported in electrical cables: it thus makes sense to measure how much charge is passing through the cable at a given point in a second. Such a measure is called {\it current}: charge per unit surface per unit time. For any portion of space one may define an {\it electric field}, consisting of a vector at each point representing the (directed) force acting on a unit charge (think of the analogy with a gravitational field, representing force on a unit mass). The {\it voltage} is the potential energy of a unit charge in the electric field. Finally, the {\it power} is measured as voltage multiplied by current. 

The transportation of power from generating plants (or {\it generators}) to users (or {\it loads}) occurs by means of a network, called {\it power grid}, the nodes of which are called {\it buses} and the links {\it lines} (or, sometimes, {\it branches}). Buses may be hubs, relays, buildings or entire corporations; lines consist of electrical cables linking two buses. Generators are assigned to buses, so that a bus may be a load, but may also, at the same time, host a varying number of generators.

Electrical networks may transport Direct Current (DC) or Alternating Current (AC). Small electronic devices, such as radios, televisions, computers, are usually driven by DC. The transportation of power on any geographical area is typically based on AC. This occurs naturally when transforming kinetic energy, e.g.~from a waterfall, to a rotating wheel which moves another wheel with some attached magnets. While the magnetic wheel rotates, the alternating magnetic field induces an alternating electrical field (in accordance with Maxwell's equations) which in turn induces an alternating current in an appropriately placed coiled cable. The frequency of the rotation goes from 50Hz to 60Hz depending on country.

While charge (in the form of electrons) actually moves along the cable in DC, in AC the charge moves very little. It ``oscillates'' in the cable, according to the alternating nature of the electrical field it is subjected to, and generates an ``energy wave'' which can be interpreted as travelling along the cables. This wave induces a power ``flow'' along the lines. The power is induced by the current and the voltage difference between the buses incident to the line. For a line $\ell$ incident to buses $b,a$, power is ``injected'' in $\ell$ at $b$ (in one direction) and at $a$ (in the other direction, with directions changing according to the field oscillations).

Electricity production stakeholders have an interest in generating sufficient power to satisfy demand. In doing so, they also have to make sure that the power can be transported from the generating bus, over the lines of the electrical networks, to the load bus. The ACOPF is supposed to achieve this purpose.

\LEO{The rest of this paper is organized as follows. We give an essential literature review in Sect.~\ref{s:litrev}. We explain how a time-dependent problem can be modelled using formulations that do not depend on time in Sect.~\ref{s:dyn2stat}. We discuss the basic modelling of the ACOPF using MP in Sect.~\ref{s:acopfmodel}. We propose some formulations in complex numbers in Sect.~\ref{s:complex}, and explain how to obtain formulations in real numbers in Sect.~\ref{s:real}. Lastly, resting on the material introduced so far, we discuss in Sect.~\ref{s:rel} some recent contributions of the theory of ACOPF relaxations, which is important in view of solving ACOPF instances of large size.}

\section{Literature review}
\label{s:litrev}

\LEO{In this section we give an essential literature review of the ACOPF and its relationship with MP. A technical literature review section about a few works specific to ACOPF relaxations is given in Sect.~\ref{s:rel}.}

\subsection{Generalities}
Early formulations of the (feasibility-only version of the) ACOPF date back to the mid-twentieth century. Power flow equations, in their complex formulation, were stated in \cite{VanNess1961}, and solved using an implementation of Newton's method. Later on, this implementation was improved in \cite{tinney1967} using sparse matrix techniques, which allowed the solution of the problem on larger networks.

A comprehensive discussion of the power flow equations is developed in Power Systems textbooks, such as \cite{andersson,bergen-book,overbye-book,Monticelli1999}. The parameters involved in these equations, including transformer, phasors and other branch elements, are discussed in detail in \cite[Ch.~4]{Monticelli1999}; this allows the formulation of power flow equations in complex and polar coordinates. The books \cite{bergen-book,overbye-book} are devoted to the physics behind the model, with technical details of power system operations.

\subsection{MP formulations}
The optimization aspect of the problem was first introduced as ``Economic Dispatch'' in \cite{Carpentier1962}, followed by a survey of the state the art of the Optimal Power Flow (OPF) problem in DC in \cite{carpentier1979}. Other extensive surveys \cite{huneault,survey1993-1,survey1993-2} of selected literature until early nineties follow the evolution of the OPF and related solution methodologies. The more recent surveys \cite{rebennack1,rebennack2} give an overview of existing MP formulations in qualitative terms, and focus on optimization methods (deterministic, non-deterministic, hybrid) for solving them.

Many algorithms relying on Interior Point Methods (IPM) have been developed in the last three decades. An IPM for nonlinear programming based on perturbing the KKT conditions of the rectangular formulation was described in \cite{wei1998}. The authors in \cite{capitanescu2007,torres1998} propose primal-dual IPM for nonlinear programs suited for the ACOPF problem, using as well the rectangular formulation. Whereas \cite{rider2004,wangpower} develop step length control techniques on the polar OPF formulation. The techniques in \cite{wangpower} were implemented in the open-source \textsc{Matpower} package \cite{matpower} for \textsc{Matlab}.

The study of OPF solutions on tree networks (also known as ``radial networks'' in the OPF literature) has allowed the development of interesting techniques. For instance, the change of variables used in the Jabr's relaxation (see Sect.~\ref{s:jabr}) was introduced in \cite{gomez1999}, where a Newton's method is proposed to solve the OPF problem in tree networks. An efficient IPM for conic quadratic programming applied to power flow equations on tree networks was later introduced in \cite{jabr} and then extended to more general networks \cite{jabr2007,jabr2008}.

A new mixed-combinatorial solution method for tree networks without transformers, based on graph reduction and expansion operators on the tree graph, was proposed in \cite{beckopf}. The reduction step repeatedly contracts leaf nodes in the star of some vertex $v$ in the tree to $v$ itself. At the same time, it updates voltage bounds at $v$ so that they are feasible w.r.t.~the ACOPF constraints imposed at the leaves. Such reductions are carried out until the only vertex left is the tree root. The opposite operation re-expands the tree while keeping voltage values feasible. The whole process appears to be similar to the well-known Feasibility-Based Bounds Tightening algorithm \cite{messine,fbbt-cocoa10} used in spatial Branch-and-Bound \cite{shectman,couenne}. 

\subsection{Relaxations}
Semidefinite Programming (SDP) relaxations based on the ACOPF formulation in rectangular coordinates have also been studied. Sufficient conditions are exhibited in \cite{lavaei} for ACOPF instances to have zero optimality gap with respect to their SDP relaxations. Moreover, optimal solutions of the original problem can be obtained from solutions of the SDP relaxation. However, \cite{lesieutre2011} proves that this SDP relaxation is not exact for a specific example of a cycle of three buses and, consequently, that it would fail on larger networks including cycles. A relaxation of the ACOPF based on Lasserre's moments hierarchy \cite{putinar2011} from the rectangular $V$-formulation of the problem is shown in \cite{molzahn2014-1,molzahn2014-2} to improve results obtained by previous relaxations such as \cite{lavaei}, where the duality gap is zero under specific conditions, the relaxations in \cite{molzahn2014-1,molzahn2014-2} improve this gap for a more general set of instances.

An extensive survey of convex formulations of the ACOPF can be found in \cite{low2013,low2014-1,low2014-2}. Moreover, \cite{coffrin2016} compares ---  theoretically and numerically --- quadratic convex relaxations derived from the complex formulation of the ACOPF problem with SDP and second order cone relaxations. A Second-Order Cone Programming (SOCP) relaxation in the $2\times 2$ minors of the Hermitian matrix variable representing voltage in rectangular coordinates is proposed in \cite{kocuk2018}, and compared with state-of-the-art SDP relaxations.

More details about ACOPF relaxations will be given in Sect.~\ref{s:rel}.

\subsection{Grid security}
The success of concrete use cases of the ACOPF critically depends on the choice of formulation. Vulnerability of the power grid, for example, has been a major concern in the past decade after serious cyber-physical attacks affecting large geographical zones \cite{sridhar2012}. The security of networks and the repercussions of line failures, such as cascades and consequential blackouts, are discussed in detail in  \cite{bienstock-acopf}. False data injection attacks \cite{chaojun2015,liu2017} aim at studying the potential threat of a cyber attack consisting in modifying the voltage measurements of the grid, which would trick the network controller into taking wrong decisions. More complicated (but possible) attacks include a physical alteration of the grid, such as disconnecting some lines or disturbing the load and generation of a zone of the grid \cite{soltan2017,bienstock2019}.

These attacks rely on the ability that a network controller has to recover or estimate the status of the system --- voltages, flows, generation and loads --- from the measurement of a (reduced) subset of these physical quantities. Since the actions that the network controller subsequently takes depend on the estimation that he or she makes of the grid status, any possible error on this process can be crucial in order to maintain the system stability. A stochastic defense mechanism that randomly perturbs the power generation in order to unmask the effect of sophisticated attacks is proposed in \cite{bienstock2019}.

\LEO{We note that many of these formulations employ binary or integer variables in order to model attacks or other vulnerabilities.}


\subsection{Network design}
\LEO{Another family of ACOPF formulation variants that depend on binary variables are those derived from network design. To the best of our knowledge, the first paper exhibiting computational results for the ACOPF with binary variables used for design purposes is \cite{panciatici-acopf}, where binary variables are used to switch generators and shunts on and off. Improved formulations with binary variables for switching generators on and off were proposed in \cite{salgado3,salgado4}. The ACOPF is \textbf{NP}-hard \cite{acopf_nphard_orl}, and remains hard even when the goal is to minimize the number of active generators \cite{salgado4}.

The rest of this survey will be concerned with formulations in continuous variables only.
}

\section{Dealing with the time dependency}
\label{s:dyn2stat}
From Sect.~\ref{s:intro}, it should be clear that optimized generation of AC power in power grids is in fact a dynamic problem, also known as a {\it mathematical control} problem. Time plays a factor insofar as energy is expressed as a wave having frequency $\omega$. Current, voltage and power values oscillate with the same frequency $\omega$ according to the corresponding relationships. We denote current on the line $\ell$ adjacent to $b,a$ by $I_{ba}$, voltage at bus $b$ by $V_b$ and power injected in $\ell$ at $b$ by $S_{ba}$. As functions of time $t$, we have the following relationships \cite[p.~3]{bienstock-acopf}:
\begin{eqnarray}
  V_b(t) &=& V^{\max}_b\cos(\omega t + \theta_b) \nonumber \label{Vt} \\
  I_{ba}(t) &=& I^{\max}_{ba}\cos(\omega t + \phi_{ba}) \nonumber \label{It} \\
  S_{ba}(t) &=& V_b(t)I_{ba}(t),\label{SVIt}
\end{eqnarray}
where $V^{\max}_b$, $I^{\max}_{ba}$ are amplitudes and $\theta_b$, $\phi_{ba}$ are phases of voltage at a bus $b$ and of current at an incident line leading to bus $a$ as functions $V_b(t)$, $I_{ba}(t)$ of time.

By Eq.~\eqref{SVIt}, using basic trigonometric relations, we obtain:
\begin{eqnarray*}
  S_{ba}(t) &=& V^{\max}_bI^{\max}_{ba}\cos(\omega t + \theta_b)\cos(\omega t + \phi_{ba}) =\\
  &=& \frac{1}{2} V^{\max}_b I^{\max}_{ba}(\cos(\theta_b-\phi_{ba}) + \cos(2\omega t + \theta_b + \phi_{ba})).
\end{eqnarray*}
Formulating a time-dependent ACOPF is certainly possible using these time-dependent quantities, but it would be infeasible to find a solution for real-life cases using current computational technology. Instead, we consider steady state average values over a period $2\pi/\omega$, which we simply indicate with $I,V,S$ without the dependence on $t$. For the usual bus $b$ and its incident line linking it to bus $a$, we then obtain \cite[Eq.~(1.3)]{bienstock-acopf}:
\begin{eqnarray}
  S_{ba} = \frac{1}{2} V_b^{\max} I_{ba}^{\max}\cos(\theta_b - \phi_{ba}),\label{SVI}
\end{eqnarray}
to which there correspond averages for $V,I$ too.

Restricting the analysis of a whole function $S(t)$ to its average $S$ causes a loss of information which is deemed excessive. An ACOPF only defined on the averages $I,V,S$ apparently fails to capture many phenomena that are important to robust power grid design. An acceptable compromise is reached by considering $I,V,S$ complex instead of real quantities. The averages computed above are then considered their real parts $\re{I},\re{V},\re{S}$; and their imaginary parts $\im{I},\im{V},\im{S}$ provide a further piece of (static) information about the dynamics of $I(t),V(t),S(t)$ as functions of time.

We therefore write cartesian and polar representations of voltage, current and power as complex quantities:\footnote{We remark that most of the power grid literature uses $i$ to indicate current, and therefore resorts to $j$ to indicate $\sqrt{-1}$. We chose to keep notation in line with mathematics and the rest of the physical sciences, namely we use $i=\sqrt{-1}$, and employ $I$ to denote current.}
\begin{eqnarray*}
  V_b &=& \re{V}_b + i\im{V}_b = \frac{V_b^{\max}}{\sqrt{2}}e^{i\theta_b} \\
  I_{ba} &=& \re{I}_{ba} + i\im{I}_{ba} = \frac{I_{ba}^{\max}}{\sqrt{2}}e^{i\phi_{ba}} \\
  S_{ba} &=& \re{S}_{ba} + i\im{S}_{ba} = |S_{ba}| e^{i\psi_{ba}},
\end{eqnarray*}
where $\psi_{ba}$ is the phase for power. We now reformulate Eq.~\eqref{SVI} as
\begin{eqnarray*}
  \re{S}_{ba} &=& |V_b||I_{ba}|\cos(\theta_b-\phi_{ba}) = \re{(V_b\conj{I_{ba}})},
\end{eqnarray*}
where $\conj{x}=\re{x}-i\im{x}$ is the {\it complex conjugate} of $x$, and $|x|=\sqrt{x\conj{x}}=\sqrt{(\re{x})^2+(\im{x})^2}$ is the {\it modulus}, of any $x\in\mathbb{C}$.

Thus, it makes sense to define an ``imaginary power''
\[ \im{S}_{ba} = \im{(V_b\conj{I_{ba}})}, \]
yielding a complex power $S_{ba}=\re{S}_{ba} + i\im{S}_{ba}$. In the power grid literature, real power is known as {\it active power}, while imaginary power is known as {\it reactive power}.

\subsection{Change of coordinates}
\label{s:coordchange}
In the rest of this paper, we will construct various ACOPF formulations based on $I_b,V_{ba},S_{ba}\in\mathbb{C}$ for any bus $b$ and any line $\ell$ incident to $b,a$. In particular, we will use both the cartesian and the polar representations of complex numbers. In this section we recall the nonlinear transformation relations between these representations.

Consider $x=\re{x}+i\im{x}\in\mathbb{C}$ expressed in cartesian coordinates. The polar representation of $x$ is $\alpha e^{i\vartheta}=\alpha\cos\vartheta+i\alpha\sin\vartheta$, where $\alpha$ is the magnitude and $\vartheta$ is known as the {\it angle} or {\it phase} of the complex number (which is itself also called {\it phasor}). 

The change of coordinates from cartesian to polar representations (and vice versa) is a nonlinear relationship between $\re{x},\im{x}$ and $\alpha,\vartheta$, as follows:  \\
\begin{equation*}
  \begin{array}{rclrcl}
  \re{x} &=& \alpha\cos\vartheta \qquad & \qquad \alpha &=& \sqrt{(\re{x})^2 + (\im{x})^2}  \\
  \im{x} &=& \alpha\sin\vartheta \qquad & \qquad \vartheta &=& \arccos(\re{x}/\alpha) = \arcsin(\im{x}/\alpha),
  \end{array}
\end{equation*}
where we take the positive sign of the square root for $\alpha$. We further remark that the following identities are often used in complex derivations:
\begin{eqnarray}
  \forall x\in\mathbb{C} \qquad 2\re{x} &=&  x+\conj{x} \label{xr} \\
  \forall x\in\mathbb{C} \ \;  -2\im{x} &=& i(x-\conj{x}). \label{xc}
\end{eqnarray}

\section{Modelling the ACOPF}
\label{s:acopfmodel}
MP formulations for the ACOPF are unlike every other formulation we have ever seen, in that it requires an unusual amount of effort to understand, and an inordinate amount of debugging in order to implement. Many OR researchers and practitioners do not use complex numbers in their normal line of work, so this is part of the difficulty. Another part is the way the input data is presented and stored, which may be natural to electrical engineers, but certainly did not seem natural to us. In this section we will do our best to explain the modelling difficulties away, and to warn the reader against the implementation pitfalls we found.

We remark that electrical engineers themselves do not all agree on the way to approximate dynamic behaviour by static quantities, nor on the notation used. We refer to the well-known, high quality, open-source and {\it de facto} standard-establishing \textsc{Matlab} software {\sc MatPower} \cite{matpower}, as well as to its user manual \cite{matpower7}, as a reference to what we mean by the expression of an ACOPF MP formulation, and its related notation, by electrical engineers. Other examples of typical notation and formulation style in use in the AC power community are given in \cite{oneill1}. 

\subsection{The power grid as a graph}
\label{s:digraph}
The power grid is a network of buses interconnected by lines. The standard abstract entity used to model networks is a graph. In this case, however, there are a few unusual modelling issues. 

\begin{enumerate}
\item The first modelling issue is that a line, which models a cable, is not always an undirected edge: whenever the line $\ell$ between buses $b,a$ has a transformer close to the $b$ end, the current and injected power on $\ell$ at $b$ is different from that at $a$ (in DC, by contrast, the current injected on $\ell$ at $b$ is equal to the negative of the current injected at $a$). This, in general, points us towards a directed graph, or {\it digraph}, where each line is modelled as a pair of anti-parallel arcs $\{(b,a),(a,b)\}$.
\item We remark that if $\ell$ hosts a transformer, $(b,a)$ and $(a,b)$ have different current and injected power values; moreover, they have voltage differences with same magnitude and opposite signs. This raises a second issue. The relationship of current and voltage for the two antiparallel arcs is expressed by a vectorial equation in two components (one for $(b,a)$ and one for $(a,b)$), where a $2$-component current vector $I_{ba}$ is given as a $2\times 2$ non-symmetric admittance matrix $\mathbf{Y}_{ba}$ multiplying a $2$-component voltage vector $\mathbf{V}_{ba}$ (voltage at $b$ and voltage at $a$):
\begin{equation}
  \mathbf{I}_{ba} = \mathbf{Y}_{ba} \mathbf{V}_{ba} \label{IYV1}
\end{equation}
(if there is a transformer on the line from $b$ to $a$, by convention it is associated with the first row of the matrix $\mathbf{Y}_{ba}$ in Eq.~\eqref{IYV1}). This modelling technique also applies to injected power because of the equation $S_{ba}=V_b\conj{I_{ba}}$. We shall also see that one of the terms in the power flow equations (see Sect.~\ref{s:informal}) is the sum $\sum_{(b,a)} S_{ba}$ of injected powers $S_{ba}$, at a bus $b$, over all the lines incident to $b$. For obvious topological reasons, this sum never involves pairs of antiparallel arcs $(b,a),(a,b)$. Instead, it may involve arcs $(b,a_1)$ (for some bus $a_1$ adjacent to $b$) where the transformer is at $b$, and arcs $(b,a_2)$ where the transformer is at $a_2$ (for some other bus $a_2$ adjacent to $b$). This means that we must pick the first component of the $2$-vector current-voltage equation Eq.~\eqref{IYV1} for $(b,a_1)$ and the second component of Eq.~\eqref{IYV1} for $(a_2,b)$, which requires considerable care with handling indices. This issue is discussed in more depth in Sect.~\ref{s:modissue}. 
\item A third issue is given by the fact that sometimes parallel cables transport power between two buses in the case when one would be insufficient for the power demand. This means there may be parallel lines $\ell_1,\ldots,\ell_p$ between two buses $b,a$, and that each must be modelled as a pair of anti-parallel arcs. This is an unusual setting insofar as MP formulations on graphs go, which contributes to the modelling difficulty associated with the ACOPF. We remark that each arc going from $b$ to $a$ has different values of power and current, but the same voltage difference associated to $b$ and $a$. Since the equations regulating power and current are nonlinear, it is not possible to model sheaves of parallel arcs by a single arc which aggregates the values of each arc in the sheaf.
\end{enumerate}

\subsection{The $\pi$-model of a line}
\label{s:linemodel}
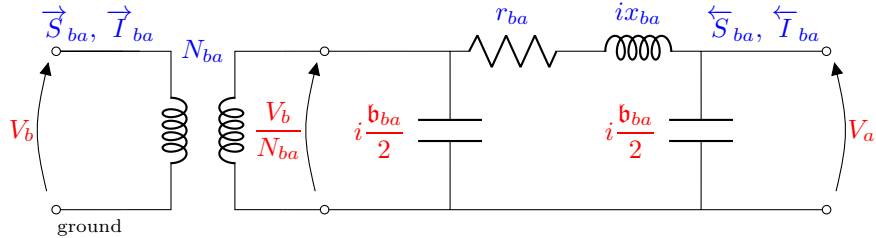
\begin{figure}[!ht]
\begin{center}
\begin{circuitikz}[scale=1.1]
\node[transformer,left=5mm] (T) {};
\draw
(T.base) node{\bluem{$N_{ba}$}} 
(T.A1) to [short, -o] ++ (-0.8,0) coordinate (in0+)
(T.A2) to [short, -o, l={\scriptsize ground}] ++ (-0.8,0) coordinate (in0-)
(T.B1) to [short, -o] ++ (0.5,0) coordinate (in+)
(T.B2) to [short, -o] ++ (0.5,0) coordinate (in-)
(in0+) to [open, v<=\redem{$V_b$}] (in0-) 
(in+) to [open, v<=\redem{$\dfrac{V_b}{N_{ba}}$}] (in-) 
(in0+) to [short, l^={\bluem{$\vec{S}_{ba}$, $\vec{I}_{ba}$}}] ++ (1,0)
(in+) to [short] ++ (1.5,0) coordinate (c1+)
to [short,C,l_=\redem{$i\dfrac{\mathfrak{b}_{ba}}{2}$}] (c1+ |- in-) 
(c1+) to [R, l=\bluem{$r_{ba}$}] ++ (1.5,0) 
to [L, l=\bluem{$ix_{ba}$}] ++ (1.5,0) coordinate (c2+) 
to [short, l=\bluem{$\cev{S}_{ba}$, $\cev{I}_{ba}$}, -o] ++ (1.5,0) coordinate (out+)
to [open, v^<=\redem{$V_{a}$}] (out+ |- in-) coordinate (out-) 
(c2+) to [C, l_=\redem{$i\dfrac{\mathfrak{b}_{ba}}{2}$}] (c2+ |- in-) 
(out-) to [short,o-] (in-); 
\end{circuitikz}
\end{center}
\caption{The $\pi$-model of a line.}
\label{f:pimodel}
\end{figure}

We introduce some of the parameter and decision variable symbols in the ACOPF by means of the so-called ``$\pi$-model'' of a line (Fig.~\ref{f:pimodel}).


The graphical representation of a line consists of two horizontal parallel segments. The above segment represents the line, going from bus $b$ (on the left) to bus $a$ (on the right). The segment below corresponds to the ground: if we were using DC, electrons flowing from $b$ to $a$ would need to come back from $a$ back to $b$ using the ground.

The pair of vertical parallel coils on the left represents a transformer installed at $b$: the parameter $N_{ba}\in\mathbb{C}$, the ``ratio'' of the AC transformation, is sometimes indicated by $N_{ba}:1$, but $N_{ba}$ is {\it not} necessarily an integer; in general, it is a complex number. It is usually expressed in its polar representation $N_{ba}=\tau_{ba}e^{i\nu_{ba}}$, where $\tau_{ba}$ is the magnitude and $\nu_{ba}$ the angle.

Next, going from left to right, we find the voltage $V_b$, expressed as a voltage difference between potential energy at $b$ and at the ground. Above the top segment we find a sequence of symbols: the injected power $\vec{S}_{ba}$ at $b$, the current $\vec{I}_{ba}$ at $b$, then the {\it series impedance} $y_{ba}=r_{ba}+ix_{ba}$ of the line $\ell$ between $b$ and $a$, and then the corresponding injections $\cev{S}_{ba}, \cev{I}_{ba}$ at $a$. Below, we find the complex terms $i\mathfrak{b}_{ba}/2$, where $\mathfrak{b}_{ba}$ is the {\it line charging susceptance} (an interaction of the line with the ground).

The $2\times 2$ matrix $\mathbf{Y}_{ba}$ referred to above is defined follows.
\begin{equation}
  \mathbf{Y}_{ba}=\left(\begin{array}{cc} Y_{ba}^{\mathsf{ff}} & Y_{ba}^{\mathsf{ft}} \\ Y_{ba}^{\mathsf{tf}} & Y_{ba}^{\mathsf{tt}} \end{array}\right)=\left(\begin{array}{cc} (\frac{1}{r_{ba}+ix_{ba}}+i\frac{\mathfrak{b}_{ba}}{2})/\tau_{ba}^2 & -\frac{1}{(r_{ba}+ix_{ba})\tau e^{-i\nu_{ba}}} \\
    -\frac{1}{(r_{ba}+ix_{ba})\tau_{ba} e^{i\nu_{ba}}} & \frac{1}{r_{ba}+ix_{ba}} + i\frac{\mathfrak{b}_{ba}}{2}\end{array}\right), \label{Yba}
\end{equation}
where the suffixes $\mathsf{ff}$, $\mathsf{ft}$, $\mathsf{tf}$, $\mathsf{tt}$ stand for ``from-from'', ``from-to'', ``to-from'', and ``to-to'', which are a reminder for the direction referring to $b,a$. Such suffixes help addressing the second issue mentioned in Sect.~\ref{s:digraph}, as we shall also see in Sect.~\ref{s:modissue} below.

We note that $\mathbf{Y}_{ba}$ is sometimes also expressed in function of $z_{ba}=\frac{1}{y_{ba}}$ as follows:
\begin{equation*}
  \mathbf{Y}_{ba}=\left(\begin{array}{cc} (z_{ba}+i\frac{\mathfrak{b}_{ba}}{2})/\tau_{ba}^2 & -z_{ba}/(\tau e^{-i\nu_{ba}}) \\
    -z_{ba}/(\tau_{ba} e^{i\nu_{ba}}) & z_{ba} + i\frac{\mathfrak{b}_{ba}}{2}\end{array}\right).
\end{equation*}

From Eq.~\eqref{IYV1} and \eqref{Yba}, we see the reason why the current on a line $\ell$ between $b$ and $a$ and hosting a transformer at $b$ is different depending on whether it flows from $b$ or $a$. We have:
\begin{eqnarray}
  \vec{I}_{ba} &=& Y_{ba}^{\mathsf{ff}}V_b + Y_{ba}^{\mathsf{ft}} V_a \label{IYV2} \\
  \cev{I}_{ba} &=& Y_{ba}^{\mathsf{tf}}V_b + Y_{ba}^{\mathsf{tt}} V_a,\label{IYV3}
\end{eqnarray}
and it is clear from Eq.~\eqref{Yba} that Eq.~\eqref{IYV2} and Eq.~\eqref{IYV3} are different. On the other hand, Eq.~\eqref{Yba} also shows us that, if no transformer is installed at the line then $N_{ab}=1$, which makes $\mathbf{Y}_{ba}$ symmetric. If, moreover, the line charging susceptance $\mathfrak{b}_{ba}$ is zero, we have $\vec{I}_{ba}=-\cev{I}_{ba}$: namely, the alternating current on the line behaves like direct current.

In this section we have denoted directed power and currents by means of arrows over the corresponding symbols $\vec{S}_{ba},\cev{S}_{ba},\vec{I}_{ba},\cev{I}_{ba}$ for clarity w.r.t.~a given line $\ell$ between $b,a$. Since our formulations below are based on a digraph representation of the network, we shall use the corresponding (directed) arc expressions $S_{ba},S_{ab},I_{ba},I_{ab}$. Moreover, if there are $p_{ba}$ parallel lines between $b$ and $a$, we use symbols $S_{bah},S_{abh},I_{bah},I_{abh}$ for $h\le p_{ba}$; the symbols $r,x,\tau,\nu,\mathfrak{b}$ are all indexed by $b,a,h$ (for $h\le p_{ba}$); and, consequently, this also holds for the symbols $y,z,Y^{\mathsf{ff}},Y^{\mathsf{ft}},Y^{\mathsf{tf}},Y^{\mathsf{tt}},\mathbf{Y}$. 

\subsection{Informal description of ACOPF formulations}
\label{s:informal}
MP formulations for the ACOPF consist of:
\begin{itemize}
\item an objective function, which usually minimizes the cost of the generated power;
\item a set of bound constraints:
  \begin{itemize}
  \item on the voltage magnitude;
  \item on the difference of voltage angles between adjacent buses;
  \item on the maximum power (or of current) magnitude injected on a line;
  \item on the generated power;
  \end{itemize}
\item a set of explicit constraints:
  \begin{itemize}
  \item on the power balance at each bus (also called ``power flow equations'');
  \item on the definition of power in function of voltage and current;
  \item on the definition of current w.r.t.~voltage (an AC generalization of Ohm's equation, see Eq.~\eqref{IYV1}).
  \end{itemize}
\end{itemize}
There may be other technical constraints, depending on the network at hand, and the application requiring the solution of the ACOPF. The {\sc MatPower} software \cite{matpower}, for example, provides a ``reference'' status for each bus. A ``reference bus'' has the voltage angle (a.k.a.~phase) set to zero.

\section{Complex formulations}
\label{s:complex}
In this section we shall discuss some ACOPF formulations where the decision variables are in $\mathbb{C}$.

\subsection{The $(S,I,V)$-formulation}
\label{s:sivform}
We first present what we consider to be the most basic ACOPF formulation, in terms of power $S$, current $I$ and voltage $V$.

\subsubsection{Sets, parameters and decision variables}
\label{s:setparamdecvar1}
Index sets and parameters represent the input of a MP formulation. The decision variables will encode the solution after an appropriate algorithm has solved the formulation.

We consider a multi-digraph $G=(B,L)$ where $B$ is the set of buses, and $L$ is the set of arcs representing the lines. We assume $|B|=n$ and $|L|=m$. $L$ is partitioned in two sets $L_0,L_1$ with $|L_0|=|L_1|$, which makes $|L|$ even: this is consistent with the fact that ever line is represented by two anti-parallel arcs. For the $h$-th line $\ell_h$ between $b,a$, for $h\le p_{ba}$, represented by a set of anti-parallel arcs $\{(b,a,h),(a,b,h)\}$, exactly one between the anti-parallel arcs, say $(b,a,h)$, belongs to $L_0$, and the other, $(a,b,h)$, belongs to $L_1$. In particular, every line $\ell_h$ with a transformer at $b$ is oriented so that $(b,a,h)\in L_0$ and $(a,b,h)\in L_1$.

We also consider a set $\mathscr{G}$ of generators, partitioned into (possibly empty) subsets $\mathscr{G}_b$ for every $b\in B$. The generators in $\mathscr{G}_b$ are those that are assigned to bus $b$.

We consider the following parameters.
\begin{enumerate}
\item The objective function is a real polynomial of degree $d$ of (the real part of) power generated at $g\in\mathscr{G}$. The coefficients $c_{g0},\ldots,c_{g,d-1}$ are given for all $g\in\mathscr{G}$.
\item Voltage.
  \begin{itemize}
  \item At each bus $b\in B$ the voltage magnitude $|V_b|$ is constrained to lie in a given real range $[\underline{V}_b,\overline{V}_b]$. 
  \item At each arc $(b,a,h)\in L_0$ the voltage phase difference $\theta_b-\theta_a$ between $b$ and $a$ is constrained to lie in given real ranges $[\underline{\eta}_{bah},\overline{\eta}_{bah}]$. Since the voltage phase difference w.r.t.~$(a,b,h)$ is simply the same for $(b,a,h)$ in absolute value, we do not need to impose these constraints for arcs in $L_1$. In fact, it suffices to impose the most restrictive bounds for each unordered pair $\{b,a\}$.
  \item A chosen bus index $r\in B$ is designated as ``reference'', which entails having voltage phase $\theta_r=0$; by setting $\re{V}_r+i\im{V}_r=|V_r|e^{i\theta_r}$, we see that $\theta_r=0$ implies $\im{V}_r=0$ and $\re{V}_r=|V_r|$.
  \end{itemize}
\item Power.
  \begin{itemize}
  \item The power demand at bus $b\in B$ is denoted $\tilde{S}_b\in\mathbb{C}$; for technical reasons, there can be buses with negative demand.
  \item At each arc $(b,a,h)\in L$ the magnitude $|S_{bah}|$  of the power injected on the line is constrained to be bounded above by a real scalar $\bar{S}_{bah}\ge 0$. This bound does not depend on the injection direction, so $\bar{S}_{bah}=\bar{S}_{abh}$.
  \item At each generator $g\in\mathscr{G}_b$ installed at bus $b\in B$ the power at $g$ is constrained to be within complex ranges $[\underline{\mathscr{S}}_{g},\overline{\mathscr{S}}_{g}]$ (meaning that the real part is in $[\re{\underline{\mathscr{S}}}_{g},\re{\overline{\mathscr{S}}}_{g}]$ and the imaginary part is in $[\im{\underline{\mathscr{S}}}_{g},\im{\overline{\mathscr{S}}}_{g}]$).    
  \end{itemize}
\item Current.
\begin{itemize}
\item For every arc $(b,a,h)\in L_0$ we are given the matrix $\mathbf{Y}_{bah}\in\mathbb{C}^{2\times 2}$. The first row of each $\mathbf{Y}_{bah}$ defines current in function of voltage for the direction from $b$ to $a$, while the second row handles the opposite direction, see Eq.~\eqref{IYV2}-\eqref{IYV3}. So we do not need to define these matrices over arcs in $L_1$.
\item At each arc $(b,a,h)\in L$ the current magnitude $|I_{bah}|$ may be constrained to be bounded above by a real scalar $\bar{I}_{bah}\ge 0$, such that $\bar{I}_{bah}=\bar{I}_{abh}$. This bound is usually enforced as an alternative to the injected power bound $\bar{S}$, see Sect.~\ref{s:vonly}.
\end{itemize}
  \item The {\it shunt admittance} at bus $b\in B$, related to an interaction with the ground, is $A_b\in\mathbb{C}$.
\end{enumerate}

We consider the following decision variables.
\begin{enumerate}[(a)]
\item The complex voltage at each bus $b\in B$ is denoted $V_b\in\mathbb{C}$.
\item The complex current at each arc $(b,a,h)\in L$ is denoted $I_{bah}\in\mathbb{C}$. For each pair of antiparallel arcs $(b,a,h)\in L_0$ and $(a,b,h)\in L_1$, we define the vector $\mathbf{I}_{bah}=\transpose{(I_{bah},I_{abh})}\in\mathbb{C}^2$.
\item The complex power at each arc $(b,a,h)\in L$ (injected at $b$) is denoted $S_{bah}\in\mathbb{C}$. For each pair of antiparallel arcs $(b,a,h)\in L_0$, $(a,b,h)\in L_1$, we define the vector $\mathbf{S}_{bah}=\transpose{(S_{bah},S_{abh})}\in\mathbb{C}^2$.
\item For a generator $g\in \mathscr{G}_b$ installed at bus $b\in B$, $\mathscr{S}_{g}$ is the complex power generated by $g$.
\end{enumerate}

\subsubsection{Objective and constraints}
The most elementary objective function employed in {\sc MatPower} \cite{matpower} is a polynomial function of degree $d$ of (real) generated power:
\begin{equation}
  \min \sum\limits_{g\in\mathscr{G}}\sum_{j=0}^d c_{g,d-j}(\re{\mathscr{S}}_g)^{d-j}.\label{obj}
\end{equation}
We also consider two simpler objective functions. One of them, with $d=2$, involves a Hermitian quadratic form (namely, with the complex square matrix $Q$ being equal to its conjugate transpose):
\begin{equation}
  \min (\hermitian{\mathscr{S}}Q\mathscr{S} + c_{1}\re{\mathscr{S}} + c_{0}\mathscr{1}), \label{obj2}
\end{equation}
where $\mathscr{S}=(\mathscr{S}_g\;|\;g\in\mathscr{G})$, $c_1=(c_{g1}\;|\;g\in\mathscr{G})$, $c_{0}=(c_{g0}\;|\;g\in\mathscr{G})$. This makes this formulation a complex Quadratically Constrained Quadratic Program (QCQP) in generated power. Another one, with $d=1$, 
\begin{equation}
  \min \sum\limits_{g\in\mathscr{G}} (c_{g1}\re{\mathscr{S}}_g + c_{g0}), \label{obj3}
\end{equation}
yields a Quadratically Constrained Program (QCP) which is linear in generated power. Eq.~\eqref{obj3} is going to be used in the voltage-only formulation of Sect.~\ref{s:vonly} so as to obtain a QCQP in voltage. 

We consider the following bound constraints.
\begin{itemize}
\item Lower/upper bounds on generated power are imposed at every generator:
  \begin{equation}
    \forall b\in B, g\in\mathscr{G}_b\quad \underline{\mathscr{S}}_{g}\le \mathscr{S}_{g}\le\overline{\mathscr{S}}_{g}. \label{genpowerbounds}
  \end{equation}
\item The upper bounds on power magnitudes (squared, since polynomial formulations are preferred) are imposed at every arc:
  \begin{eqnarray}
    \forall (b,a,h)\in L\quad |S_{bah}|^2 &\le& \bar{S}_{bah}^2, \label{powerbound}
  \end{eqnarray}  
  where we recall that $\bar{S}_{bah}=\bar{S}_{abh}$.
\item The bounds $[\underline{\eta}_{bah},\overline{\eta}_{bah}]$ on the phase difference cannot be imposed directly, as the voltage phase does not appear as a decision variable. From the polar representation $x = |x|e^{i\vartheta}$ one could argue $\vartheta= -i\ln(x/|x|)$, but this would prevent the formulation from being quadratic in its variables. Instead, we proceed as follows. We select an appropriate monotonically increasing function and write the original constraints:
  \[ \forall (b,a,h)\in L_0 \quad \underline{\eta}_{bah}\le\theta_b-\theta_a\le\overline{\eta}_{bah} \] as
  \[ \forall (b,a,h)\in L_0 \quad \tan(\underline{\eta}_{bah})\le\tan(\theta_b-\theta_a)\le\tan(\overline{\eta}_{bah}), \]
  assuming that $-\pi/4<\underline{\eta}_{bah}\le\overline{\eta}_{bah}<\pi/4$ (this assumption is justified in the following sense: either these constraints are inactive, in which case they need not be enforced, or else, in practice, $\overline{\eta}_{bah}-\underline{\eta}_{bah}$ is usually smaller than $\pi/2$).
  Next, we note that
  \begin{eqnarray*}
    \tan(\theta_b-\theta_a) &=& \frac{\sin(\theta_b-\theta_a)}{\cos(\theta_b-\theta_a)}=\frac{|V_b|\,|V_a|\sin(\theta_b-\theta_a)}{|V_b|\,|V_a|\cos(\theta_b-\theta_a)} \\
    &=& \frac{|V_b|\sin\theta_b|V_a|\cos\theta_a-|V_b|\cos\theta_b|V_a|\sin\theta_a}{|V_b|\cos\theta_b|V_a|\cos\theta_a+|V_b|\sin\theta_b|V_a|\sin\theta_a} \\
    &=& \frac{\im{V}_b\re{V}_a-\re{V}_b\im{V}_a}{\re{V}_b\re{V}_a+\im{V}_b\im{V}_a} = \frac{\im{(V_b\conj{V_a})}}{\re{(V_b\conj{V_a})}},
  \end{eqnarray*}
  whence the desired constraints can be written as:
  \begin{eqnarray}
    \forall (b,a,h)\in L_0\ \tan(\underline{\eta}_{bah}) \le& \frac{\im{(V_b\conj{V_a})}}{\re{(V_b\conj{V_a})}} & \le\tan(\overline{\eta}_{bah}) \nonumber \label{phasediffbound1} \\
    \Rightarrow\ \tan(\underline{\eta}_{bah})\re{(V_b\conj{V_a})} \le& \im{(V_b\conj{V_a})} & \le\tan(\overline{\eta}_{bah})\re{(V_b\conj{V_a})} \label{phasediffbound2}    
  \end{eqnarray}
  provided the additional constraints:
  \begin{equation}
    \forall (b,a,1)\in L_0 \quad \re{(V_b\conj{V_a})}\ge 0
    \label{phasediffboundaux}
  \end{equation}
  hold. Note that Eq.~\eqref{phasediffbound2}-\eqref{phasediffboundaux} are both quadratic in voltage. Note also that imposing the bounds on the tangent rather than on the angles requires $[\underline{\eta}_{bah},\overline{\eta}_{bah}]\subset[-\pi/2,\pi/2]$; this restriction is not problematic in practice, since anything larger really corresponds to the absence of these bounds. Note that Eq.~\eqref{phasediffboundaux} can be strengthened by requiring
  \begin{equation*}
    \forall (b,a,1)\in L_0 \quad \overline{V}_b\overline{V}_a\max\{0,\cos(\underline{\eta}_{bah}),\cos(\overline{\eta}_{bah})\}\le\re{(V_b\conj{V_a})}\le   \overline{V}_b\overline{V}_a.
  \end{equation*}
\item Lower/upper bounds on the voltage magnitude are imposed at each bus:
  \begin{equation}
    \forall b\in B \quad \underline{V}_b^2 \le |V_b|^2 \le \overline{V}_b^2.\label{voltagebound}
  \end{equation}
\item At the reference bus $r\in B$, we have $\im{V}_r=0$; as mentioned above, this makes $\re{V}_r=|V_r|$, implying
  \begin{equation}
    \im{V}_r=0\quad\land\quad \re{V}_r\ge 0. \label{reference}
  \end{equation}
\end{itemize}

We consider the following explicit constraints.
\begin{itemize}
\item The power flow equations state that, at each bus $b\in B$, the sum of complex injected powers $S_{bah}$ at bus $b$, plus the power demand $\tilde{S}_b$ at $b$, is equal to the power generated by any generators installed at $b$, plus the shunt admittance term:
  \begin{equation}
    \forall b\in B\quad \sum\limits_{(b,a,h)\in L} S_{bah} + \tilde{S}_b = -\conj{A_b}|V_b|^2 + \sum\limits_{g\in\mathscr{G}_b} \mathscr{S}_g. \label{powerflow}
  \end{equation}
  As regards shunt admittance, we remark that $-\conj{A_b}=-\re{A}_b+i\im{A}_b$. The shunt admittance $-\conj{A_b}|V_b|^2$ arises when considering the line equations in \cite[Eq.~(1.10)-(1.11)]{bienstock-acopf} applied to a fictitious line between a bus representing the ground (with associated zero voltage magnitude) and the bus $b$. In particular, it can be derived from \cite[Eq.~(1.12)]{bienstock-acopf} by setting the voltage at the ``to'' node (called $V_m$ in the cited equation) to zero.
\item The definition of power in terms of current is:
  \begin{eqnarray}
    \forall (b,a,h)\in L \quad S_{bah} &=& V_b\conj{I_{bah}}.\label{powercurrent}
  \end{eqnarray}
\item The generalized Ohm's law, which relates current to voltage, is as follows:
  \begin{eqnarray}
    \forall (b,a,h)\in L_0\quad I_{bah} &=& Y^{\mathsf{ff}}_{bah}V_b + Y^{\mathsf{ft}}_{bah}V_a\label{ohmlaw1} \\
    \forall (b,a,h)\in L_0\quad I_{abh} &=& Y^{\mathsf{tf}}_{bah}V_b + Y^{\mathsf{tt}}_{bah}V_a.\label{ohmlaw2}
  \end{eqnarray}
\end{itemize}

\subsection{A modelling issue with the power flow equations}
\label{s:modissue}
We can now provide a deeper analysis of the second modelling issue presented in Sect.~\ref{s:digraph}. The first term of Eq.~\eqref{powerflow} is the sum
\[\sum\limits_{(b,a,h)\in L} S_{bah}, \]
where $b$ is fixed by the quantifier $\forall b\in B$, and the sum ranges over all $a,h$ such that $(b,a,h)\in L$. Since $L$ is partitioned into $L_0$ and $L_1$ as explained in Sect.~\ref{s:setparamdecvar1}, some indices of terms in the sum might be arcs $(b,a,h)\in L_0$ and some others might be in $L_1$. Moreover, no pair of antiparallel arc can ever appear in the same sum. Hence, terms $S_{bah}$ with indices in $L_0$ correspond to currents $I_{bah}$ that transform to voltage by means of the first row $(Y^{\mathsf{ff}}_{ba},Y^{\mathsf{ft}}_{ba})$ of $\mathbf{Y}_{ba}$, while terms $S_{bah}$ with indices in $L_1$ correspond to currents $I_{bah}$ that transform to voltage by means of the second row $(Y^{\mathsf{tf}}_{ab},Y^{\mathsf{tt}}_{ab})$ of $\mathbf{Y}_{ab}$.

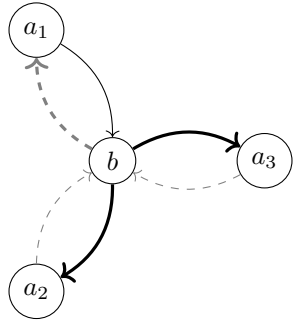
\begin{figure}[!ht]
  \begin{center}
    \begin{tikzpicture}
      \node[circle, draw=black] at (360:0mm) (center) {$b$};
      \node[circle, draw=black] at (120:2cm) (n1) {$a_1$};
      \draw[->,very thick,gray,dashed] (center) to [bend left] (n1);
      \draw[->] (n1) to [bend left] (center);
      \node[circle, draw=black] at (240:2cm) (n2) {$a_2$};
      \draw[->,very thick] (center) to [bend left] (n2);
      \draw[->,gray,dashed] (n2) to [bend left] (center);
      \node[circle, draw=black] at (360:2cm) (n3) {$a_3$};
      \draw[->,very thick] (center) to [bend left] (n3);
      \draw[->,gray,dashed] (n3) to [bend left] (center);
    \end{tikzpicture}
  \end{center}
  \caption{Example showing how the rows of $\mathbf{Y}$ matrices apply to terms in the sum of the power flow equations.}
  \label{f:pflow}
\end{figure}

The example in Fig.~\ref{f:pflow} shows a bus $b\in B$ with an adjacent neighbourhood $\{a_1,a_2,a_3\}\subset B$. Each line is represented by two anti-parallel arcs; there are no parallel lines, so we can dispense with the index $h$. 
We assume that $(a_1,b), (b,a_2), (b,a_3)$, drawn in black in Fig.~\ref{f:pflow}, are in $L_0$; the anti-parallel arcs $(b,a_1),(a_2,b),(a_3,b)$, drawn in dashed grey, are in $L_1$. The sum above is centered in $b$, so it works out to $S_{ba_1}+S_{ba_2}+S_{ba_3}$; the terms in the sum are indexed by thick arcs in Fig.~\ref{f:pflow}. The corresponding constraints in Eq.~\eqref{powercurrent} are:
\begin{eqnarray*}
  S_{ba_1} &=& V_b\conj{I_{ba_1}} \\
  S_{ba_2} &=& V_b\conj{I_{ba_2}} \\
  S_{ba_3} &=& V_b\conj{I_{ba_3}}.
\end{eqnarray*}
The corresponding constraints in Eq.~\eqref{ohmlaw1}-\eqref{ohmlaw2} are:
\begin{eqnarray*}
  I_{ba_1} &=& Y_{a_1b}^{\mathsf{tf}} V_{a_1} + Y_{a_1b}^{\mathsf{tt}} V_b \\
  I_{ba_2} &=& Y_{ba_2}^{\mathsf{ff}} V_b + Y_{ba_2}^{\mathsf{ft}} V_{a_2} \\
  I_{ba_3} &=& Y_{ba_3}^{\mathsf{ff}} V_b + Y_{ba_3}^{\mathsf{ft}} V_{a_3}.
\end{eqnarray*}
In particular, note the asymmetry between the definition of $I_{ba_1}$ and those of $I_{ba_2},I_{ba_3}$. This asymmetry is not clearly apparent in the variable indexing in the power flow equations Eq.~\eqref{powerflow}, so it may lead to modelling mistakes.

This modelling issue can be addressed by separating the sum into two sums, one over $L_0$ and the other over $L_1$, so that the first sum behaves according to Eq.~\eqref{ohmlaw1} while the second to Eq.~\eqref{ohmlaw2} reindexed so that the bus close to the transformer is called $a$:
\begin{equation}
  \forall (b,a,h)\in L_1 \quad I_{bah} = Y_{abh}^{\mathsf{tf}} V_a + Y_{abh}^{\mathsf{tt}} V_b. \label{ohmlaw2a}
\end{equation}
We remark that the values in the matrix $\mathbf{Y}_{bah}$ in Eq.~\eqref{ohmlaw2} are the very same values in the matrix $\mathbf{Y}_{abh}$ in Eq.~\eqref{ohmlaw2a}. Finally, the relevant reformulation of Eq.~\eqref{powerflow} is:
\begin{equation}
  \forall b\in B\quad \sum\limits_{(b,a,h)\in L_0} S_{bah} + \sum\limits_{(b,a,h)\in L_1} S_{bah} + \tilde{S}_b = -\conj{A_b}|V_b|^2 + \sum\limits_{g\in\mathscr{G}_b} \mathscr{S}_g. \label{powerflow2}
\end{equation}

\subsection{Voltage-only formulation}
\label{s:vonly}
This formulation is obtained from the $(S,I,V)$-formulation in Eq.~\eqref{obj3}-\eqref{ohmlaw2} by replacing injected power variables $S$ by voltage and current variables using Eq.~\eqref{powercurrent}, and then all current variables $I$ using Eq.~\eqref{ohmlaw1}-\eqref{ohmlaw2}. Instead of the power flow equations in form Eq.~\eqref{powerflow} we use Eq.~\eqref{powerflow2} according to the discussion in Sect.~\ref{s:modissue}. Although this formulation is known as ``voltage-only'', it has two types of variables: voltage $V$ and generated power $\mathscr{S}$. In all of this section, it is important to recall that $\bar{S}_{bah}=\bar{S}_{abh}$.

By simple substitution using Eq.~\eqref{ohmlaw1} and \eqref{ohmlaw2a} followed by Eq.~\eqref{powercurrent}, and recalling that $\forall x,y\in\mathbb{C}\;\conj{(x+y)}=\conj{x}+\conj{y}$ and $\conj{(xy)}=\conj{x}\conj{y}$, we have:
\begin{eqnarray}
  \forall (b,a,h)\in L_0 \quad S_{bah} &=& \conj{Y_{bah}^{\mathsf{ff}}} |V_b|^2 + \conj{Y_{bah}^{\mathsf{ft}}} V_b\conj{V_a} \label{srepl1} \\
  \forall (b,a,h)\in L_1 \quad S_{bah} &=& \conj{Y_{abh}^{\mathsf{tf}}} V_b\conj{V_a} + \conj{Y_{abh}^{\mathsf{tt}}} |V_b|^2. \label{srepl2}
\end{eqnarray}

Carrying out the replacements Eq.~\eqref{srepl1}-\eqref{srepl2} for $S$ in terms of $V$ in Eq.~\eqref{obj3}-\eqref{powerflow} yields the MP formulation below:
{
\begin{equation}
  \left.
  \begin{array}{rrcl} 
    \min & \sum\limits_{g\in\mathscr{G}} (c_{g1}\re{\mathscr{S}}_g + c_{g0}) && \\
    \forall b\in B, g\in\mathscr{G}_b & \underline{\mathscr{S}}_{g} \le \mathscr{S}_{g} &\le& \overline{\mathscr{S}}_{g} \\
    \forall (b,a,h)\in L_0 & |\conj{Y_{bah}^{\mathsf{ff}}} |V_b|^2 + \conj{Y_{bah}^{\mathsf{ft}}} V_b\conj{V_a}| &\le& \bar{S}_{bah} \\
    \forall (b,a,h)\in L_1 & |\conj{Y_{abh}^{\mathsf{tf}}} V_b\conj{V_a} + \conj{Y_{abh}^{\mathsf{tt}}} |V_b|^2| &\le& \bar{S}_{bah} \\
    \forall (b,a,h)\in L_0 & [\tan(\underline{\eta}_{bah}), \tan(\overline{\eta}_{bah})]\re{(V_b\conj{V_a})} & \ni & \im{(V_b\conj{V_a})} \\
    \forall (b,a,1)\in L_0 & \re{(V_b\conj{V_a})} &\ge& 0 \\
    \forall b\in B & \underline{V}_b^2\ \le\ |V_b|^2\ &\le& \overline{V}_b^2 \\
    & \im{V}_r = 0 \quad \land\quad \re{V}_r &\ge & 0 \\
    \forall b\in B &\sum\limits_{(b,a,h)\in L_0} (\conj{Y_{bah}^{\mathsf{ff}}} |V_b|^2 + \conj{Y_{bah}^{\mathsf{ft}}} V_b\conj{V_a}) &+& \\  &+ \sum\limits_{(b,a,h)\in L_1}(\conj{Y_{abh}^{\mathsf{tf}}} V_b\conj{V_a} + \conj{Y_{abh}^{\mathsf{tt}}} |V_b|^2) &+& \\ &+ \tilde{S}_b = -\conj{A_b}|V_b|^2 + \sum\limits_{g\in\mathscr{G}_b} \mathscr{S}_g.  && 
  \end{array}
  \right\} \label{vonly}
\end{equation}
}
The constraints of Eq.~\eqref{vonly} follow the order given in Eq.~\eqref{obj3}-\eqref{powerflow}: generated power bounds, upper bounds on injected power, phase difference bounds, voltage magnitude bounds, reference bus bounds, and power flow equations. 

Eq.~\eqref{vonly} is not quite a complex QCQP: the injected power bounds Eq.~\eqref{powerbound} either yield quartic polynomials in voltage, or the modulus of power, involving a square root. In order to obtain a complex QCQP we need to replace these constraints with corresponding constraints on the current magnitude: we replace $|S_{bah}|\le \bar{S}_{bah}$ with $|I_{bah}|^2\le\bar{I}^2_{bah}$ over all $(b,a,h)\in L_0$ (and correspondingly for $L_1$). After replacing $I_{bah}$ with the corresponding voltage terms according to Ohm's law, we obtain the following:
\begin{eqnarray}
    \forall (b,a,h)\in L_0 \quad |Y_{bah}^{\mathsf{ff}} V_b + Y_{bah}^{\mathsf{ft}} V_a|^2 &\le& \bar{I}^2_{bah} \label{Vcurrentbound1} \\
    \forall (b,a,h)\in L_1 \quad |Y_{abh}^{\mathsf{tf}} V_a + Y_{abh}^{\mathsf{tt}} V_b|^2 &\le& \bar{I}^2_{bah}. \label{Vcurrentbound2}
\end{eqnarray}
Note that Eq.~\eqref{Vcurrentbound1}-\eqref{Vcurrentbound2} are quadratic in $V$, as required.

Technically speaking, there is no loss of information or precision in replacing power magnitude bounds with current magnitude bounds, since upper bounds to injected power are often derived from upper bounds to current on the same line. On the other hand, the actual data for $\bar{I}$ may not be given. In this case we can form the relaxation given by replacing the injected power bound constraints in Eq.~\eqref{vonly} with:
\begin{eqnarray}
    \forall (b,a,h)\in L_0 \quad |Y_{bah}^{\mathsf{ff}} V_b + Y_{bah}^{\mathsf{ft}} V_a|^2 &\le& \bar{S}^2_{bah}/\underline{V}_b^2 \label{Vcurrentbound1rel} \\
    \forall (b,a,h)\in L_1 \quad |Y_{abh}^{\mathsf{tf}} V_a + Y_{abh}^{\mathsf{tt}} V_b|^2 &\le& \bar{S}^2_{bah}/\underline{V}_b^2. \label{Vcurrentbound2rel}
\end{eqnarray}
Eq.~\eqref{Vcurrentbound1rel}-\eqref{Vcurrentbound2rel} are valid constraints for Eq.~\eqref{vonly} because:
\begin{eqnarray*}
  \underline{V}_b\le |V_b| &\Rightarrow& \underline{V}_b|I_{bah}|\le|V_b|\,|I_{bah}|\Rightarrow\underline{V}_b|I_{bah}|\le\\
  &\le& |V_b|\,|\conj{I_{bah}}|=|V_b\conj{I_{bah}}|=|S_{bah}|\le\bar{S}_{bah}
\end{eqnarray*}
for all $(b,a,h)\in L_0$, whence $|I_{bah}|\le\bar{S}_{bah}/\underline{V}_b$. The argument for $(b,a,h)\in L_1$ is similar (note that $|xy|=|x|\,|y|$ for $x,y\in\mathbb{C}$). Hence Eq.~\eqref{Vcurrentbound1rel}-\eqref{Vcurrentbound2rel} provide a relaxation of Eq.~\eqref{vonly}, as claimed. Note that Eq.~\eqref{Vcurrentbound1rel}-\eqref{Vcurrentbound2rel} are quadratic in $V$, so this relaxation is also a complex QCQP.

\subsection{SDP relaxation}
\label{s:sdprel}
We derive a (complex) SDP relaxation from the complex QCQP arising from Eq.~\eqref{vonly}, and using the current magnitude bound constraints Eq.~\eqref{Vcurrentbound1}-\eqref{Vcurrentbound2}. For each $(b,a,h)\in L_0$ we have:
{
\begin{eqnarray*}
  && |Y_{bah}^{\mathsf{ff}} V_b + Y_{bah}^{\mathsf{ft}} V_a|^2 = \\
  &=& (Y_{bah}^{\mathsf{ff}} V_b + Y_{bah}^{\mathsf{ft}} V_a)\conj{(Y_{bah}^{\mathsf{ff}} V_b + Y_{bah}^{\mathsf{ft}} V_a)} \\
  &=& (Y_{bah}^{\mathsf{ff}} V_b + Y_{bah}^{\mathsf{ft}} V_a)(\conj{Y_{bah}^{\mathsf{ff}}}\conj{V_b} + \conj{Y_{bah}^{\mathsf{ft}}}\conj{V_a}) \\
  &=& |Y_{bah}^{\mathsf{ff}}|^2|V_b|^2 + |Y_{bah}^{\mathsf{ft}}|^2|V_a|^2 + Y_{bah}^{\mathsf{ff}} \conj{Y_{bah}^{\mathsf{ft}}} V_b\conj{V_a} + \conj{Y_{bah}^{\mathsf{ff}}} Y_{bah}^{\mathsf{ft}} \conj{V_b} V_a \\
  &=& |Y_{bah}^{\mathsf{ff}}|^2|V_b|^2 + |Y_{bah}^{\mathsf{ft}}|^2|V_a|^2 + Y_{bah}^{\mathsf{ff}} \conj{Y_{bah}^{\mathsf{ft}}} V_b\conj{V_a} + \conj{(Y_{bah}^{\mathsf{ff}} \conj{Y_{bah}^{\mathsf{ft}}})} \conj{(V_b\conj{V_a})} \\
  &=& |Y_{bah}^{\mathsf{ff}}|^2|V_b|^2 + |Y_{bah}^{\mathsf{ft}}|^2|V_a|^2 + Y_{bah}^{\mathsf{ff}} \conj{Y_{bah}^{\mathsf{ft}}} V_b\conj{V_a} + \conj{(Y_{bah}^{\mathsf{ff}} \conj{Y_{bah}^{\mathsf{ft}}} V_b\conj{V_a})} \\
  &=& |Y_{bah}^{\mathsf{ff}}|^2|V_b|^2 + |Y_{bah}^{\mathsf{ft}}|^2|V_a|^2 + 2\re{(Y_{bah}^{\mathsf{ff}} \conj{Y_{bah}^{\mathsf{ft}}} V_b \conj{V_a})}.
\end{eqnarray*}
}
Similarly, for each $(b,a,h)\in L_1$ we have:
{
\begin{eqnarray*}
  && |Y_{abh}^{\mathsf{tf}} V_a + Y_{abh}^{\mathsf{tt}} V_b|^2 =\\
  &=& |Y_{abh}^{\mathsf{tf}}|^2|V_a|^2 + |Y_{abh}^{\mathsf{tt}}|^2|V_b|^2 + (Y_{abh}^{\mathsf{tf}} \conj{Y_{abh}^{\mathsf{tt}}} V_a \conj{V_b}) + (\conj{Y_{abh}^{\mathsf{tf}}} Y_{abh}^{\mathsf{tt}} \conj{V_a} V_b)\\
  &=& |Y_{abh}^{\mathsf{tf}}|^2|V_a|^2 + |Y_{abh}^{\mathsf{tt}}|^2|V_b|^2 + 2\re{(Y_{abh}^{\mathsf{tf}} \conj{Y_{abh}^{\mathsf{tt}}} V_a \conj{V_b})}.  
\end{eqnarray*}
}%
We then rewrite Eq.~\eqref{vonly} as a complex QCQP as follows, where the terms for the squared current modulus have been further modified via Eq.~\eqref{xr}-\eqref{xc}. 
{
\begin{equation}
  \left.
  \begin{array}{rrcl} 
    \min & \sum\limits_{g\in\mathscr{G}} (c_{g1}\re{\mathscr{S}}_g + c_{g0}) && \\
    \forall b\in B, g\in\mathscr{G}_b & \underline{\mathscr{S}}_{g} \le \mathscr{S}_{g} &\le& \overline{\mathscr{S}}_{g} \\ [0.4em]
    \forall (b,a,h)\in L_0 & |Y_{bah}^{\mathsf{ff}}|^2|V_b|^2 + |Y_{bah}^{\mathsf{ft}}|^2|V_a|^2 &+& \\ [0.2em] &+ Y_{bah}^{\mathsf{ff}} \conj{Y_{bah}^{\mathsf{ft}}} V_b\conj{V_a} + \conj{Y_{bah}^{\mathsf{ff}}} Y_{bah}^{\mathsf{ft}} V_a\conj{V_b} &\le& \bar{I}_{bah}^2 \\ [0.4em]
    \forall (b,a,h)\in L_1 & |Y_{abh}^{\mathsf{tf}}|^2 |V_a|^2 + |Y_{abh}^{\mathsf{tt}}|^2 |V_b|^2 &+& \\ [0.2em] &+ Y_{abh}^{\mathsf{tf}} \conj{Y_{abh}^{\mathsf{tt}}} V_a \conj{V_b} + \conj{Y_{abh}^{\mathsf{tf}}} Y_{abh}^{\mathsf{tt}} V_b\conj{V_a} &\le& \bar{I}_{bah}^2 \\ [0.4em]
    \forall (b,a,h)\in L_0 & (\tan(\underline{\eta}_{bah})+i)V_b\conj{V_a} + (\tan(\underline{\eta}_{bah})-i) V_a\conj{V_b} &\le & 0 \\ [0.4em]
    \forall (b,a,h)\in L_0 & (\tan(\overline{\eta}_{bah})+i)V_b\conj{V_a} + (\tan(\overline{\eta}_{bah})-i) V_a\conj{V_b} &\ge& 0 \\ [0.4em]
    \forall (b,a,1)\in L_0 & (V_b\conj{V_a} + V_a\conj{V_b}) &\ge& 0 \\ [0.4em]
    \forall b\in B & \underline{V}_b^2\ \le\ |V_b|^2\ &\le& \overline{V}_b^2 \\ [0.1em]
    & V_r-\conj{V_r} = 0 \quad \land\quad V_r+\conj{V_r} &\ge & 0 \\ [0.4em]
    \forall b\in B &\sum\limits_{(b,a,h)\in L_0} (\conj{Y_{bah}^{\mathsf{ff}}} |V_b|^2 + \conj{Y_{bah}^{\mathsf{ft}}} V_b\conj{V_a}) &+& \\  &+ \sum\limits_{(b,a,h)\in L_1}(\conj{Y_{abh}^{\mathsf{tf}}} V_b\conj{V_a} + \conj{Y_{abh}^{\mathsf{tt}}} |V_b|^2) &+& \\ &+ \tilde{S}_b = -\conj{A_b}|V_b|^2 + \sum\limits_{g\in\mathscr{G}_b} \mathscr{S}_g,  && 
  \end{array}
  \right\} \label{vonly2}
\end{equation}
}%
so that the voltage appears linearly only as $V_r$, and quadratically as $|V_b|^2$, $V_b\conj{V_a}$, and $\conj{(V_b\conj{V_a})}=V_a\conj{V_b}$. We remark that the phase difference bound inequalities in Eq.~\eqref{vonly2} (involving tangents) follow from Eq.~\eqref{xr}-\eqref{xc}.

We now form a matrix of $n^2$ decision variables products, which we linearize by using a Hermitian matrix $X$:
{
\begin{equation*}
  V\hermitian{V} = \left(\begin{array}{cccc}
    |V_1|^2 & V_1\conj{V_2} & \cdots & V_1\conj{V_n} \\
    V_2\conj{V_1} & |V_2|^2 & \cdots & V_2\conj{V_n} \\
    \vdots & \vdots & \ddots & \vdots \\
    V_n\conj{V_1} & V_n\conj{V_2} & \cdots & |V_n|^2 
  \end{array}\right) =
  \left(\begin{array}{cccc}
    X_{11} & X_{12} & \cdots & X_{1n} \\
    X_{21} & X_{22} & \cdots & X_{2n} \\
    \vdots & \vdots & \ddots & \vdots \\
    X_{n1} & X_{n2} & \cdots & X_{nn}
  \end{array}\right)=X.
\end{equation*}
}

This directly leads to a complex SDP relaxation of the ACOPF, which consists in relaxing $X=V\hermitian{V}$ to $X\succeq_{\mathsf{H}} V\hermitian{V} \ (\star)$, \LEO{where $\succeq_{\mathsf{H}}$ denotes the Loewner order on Hermitian matrices: two Hermitian matrices $A,B$ are in the order $A\succeq_H B$ if $A-B$ is positive semidefinite (PSD). This order naturally restricts to the familiar order $A\succeq B$ to mean $A-B$ is PSD on real matrices $A,B$. The formulation is as follows:} 
{
\begin{equation}
  \left.
  \begin{array}{rrcl} 
    \min & \sum\limits_{g\in\mathscr{G}} (c_{g1}\re{\mathscr{S}}_g + c_{g0}) && \\
    \forall b\in B, g\in\mathscr{G}_b & \underline{\mathscr{S}}_{g} \le \mathscr{S}_{g} &\le& \overline{\mathscr{S}}_{g} \\ [0.4em]
    \forall (b,a,h)\in L_0 & |Y_{bah}^{\mathsf{ff}}|^2X_{bb} + |Y_{bah}^{\mathsf{ft}}|^2X_{aa} &+& \\ [0.2em] &+ Y_{bah}^{\mathsf{ff}} \conj{Y_{bah}^{\mathsf{ft}}} X_{ba} + \conj{Y_{bah}^{\mathsf{ff}}} Y_{bah}^{\mathsf{ft}} X_{ab} &\le& \bar{I}_{bah}^2 \\ [0.4em]
    \forall (b,a,h)\in L_1 & |Y_{abh}^{\mathsf{tf}}|^2 X_{aa} + |Y_{abh}^{\mathsf{tt}}|^2 X_{bb} &+& \\ [0.2em] &  + Y_{abh}^{\mathsf{tf}} \conj{Y_{abh}^{\mathsf{tt}}} X_{ab} + \conj{Y_{abh}^{\mathsf{tf}}} Y_{abh}^{\mathsf{tt}} X_{ba} &\le& \bar{I}_{bah}^2 \\ [0.4em]
    \forall (b,a,h)\in L_0 & (\tan(\underline{\eta}_{bah})+i)X_{ba} + (\tan(\underline{\eta}_{bah})-i) X_{ab} &\le & 0 \\ [0.4em]
    \forall (b,a,h)\in L_0 & (\tan(\overline{\eta}_{bah})+i) X_{ba} + (\tan(\overline{\eta}_{bah})-i) X_{ab} &\ge& 0 \\ [0.4em]
    \forall (b,a,1)\in L_0 & (X_{ba} + X_{ab}) &\ge& 0 \\ [0.4em]
    \forall b\in B & \underline{V}_b^2\ \le\ X_{bb} \ &\le& \overline{V}_b^2 \\ [0.1em]
    & V_r-\conj{V_r} = 0 \quad \land\quad V_r+\conj{V_r} &\ge & 0 \\ [0.4em]
    \forall b\in B &\sum\limits_{(b,a,h)\in L_0} (\conj{Y_{bah}^{\mathsf{ff}}} X_{bb} + \conj{Y_{bah}^{\mathsf{ft}}} X_{ba}) &+& \\  &+ \sum\limits_{(b,a,h)\in L_1}(\conj{Y_{abh}^{\mathsf{tf}}} X_{ba} + \conj{Y_{abh}^{\mathsf{tt}}} X_{bb}) &+& \\ &+ \tilde{S}_b = -\conj{A_b}X_{bb} + \sum\limits_{g\in\mathscr{G}_b} \mathscr{S}_g  && \\
    & \left(\begin{array}{cc} 1 & \hermitian{V} \\ V & X\end{array}\right) &\succeq_{\mathsf{H}} & 0.
  \end{array}
  \right\} \label{vsdp}
\end{equation}
}
The last constraint in Eq.~\eqref{vsdp} is derived from ($\star$) using the Schur complement.

A complex SDP relaxation (widely used in the literature) is obtained by relaxing the reference bus constraint and bounding the injected power on lines instead of the current on lines. \LEO{The decision variables of this formulation are the matrix $X$ and the} flow variables $S_{bah}$ (for $(b,a,h)\in L$). The formulation is as follows:
{
\begin{equation}
  \left.
  \begin{array}{rrcl} 
    \min & \sum\limits_{g\in\mathscr{G}} (c_{g1}\re{\mathscr{S}}_g + c_{g0}) && \\
    \forall b\in B, g\in\mathscr{G}_b & \underline{\mathscr{S}}_{g} \le \mathscr{S}_{g} &\le& \overline{\mathscr{S}}_{g} \\ [0.4em]
    \forall (b,a,h)\in L_0 & \conj{Y_{bah}^{\mathsf{ff}}}X_{bb} + \conj{Y_{bah}^{\mathsf{ft}}}X_{ba} &=& S_{bah}\\ [0.2em] 
    \forall (b,a,h)\in L_1 & \conj{Y_{abh}^{\mathsf{tf}}}X_{ba} + \conj{Y_{abh}^{\mathsf{tt}}}X_{bb} &=& S_{bah}\\ [0.2em]   
    \forall (b,a,h)\in L & S_{bah}\conj{S_{bah}} &\le& \overline{S}_{bah}^2\\ [0.2em]   
    \forall (b,a,h)\in L_0 & (\tan(\underline{\eta}_{bah})+i)X_{ba} + (\tan(\underline{\eta}_{bah})-i) X_{ab} &\le & 0 \\ [0.4em]
    \forall (b,a,h)\in L_0 & (\tan(\overline{\eta}_{bah})+i) X_{ba} + (\tan(\overline{\eta}_{bah})-i) X_{ab} &\ge& 0 \\ [0.4em]
    \forall (b,a,1)\in L_0 & (X_{ba} + X_{ab}) &\ge& 0 \\ [0.4em]
    \forall b\in B & \underline{V}_b^2\ \le\ X_{bb} \ &\le& \overline{V}_b^2 \\ [0.1em]
    \forall b\in B &\sum\limits_{(b,a,h)\in L} S_{bah} + \tilde{S}_b + \conj{A_b}X_{bb} - \sum\limits_{g\in\mathscr{G}_b} \mathscr{S}_g  &=& 0 \\
    & X &\succeq_{\mathsf{H}} & 0.
  \end{array}
  \right\} \label{Xsdp}
\end{equation}
}
We note that the PSD constraint $X\succeq_{\mathsf{H}}0$ is equivalent to ($\star$) for Eq.~\eqref{Xsdp} because the $V$ variables do not appear therein.

\section{Real formulations}
\label{s:real}
In this section we shall explain how to obtain real number formulations from the formulations in complex numbers discussed in Sect.~\ref{s:complex}. 

\subsection{Cartesian $(S,I,V)$-formulation}
\label{s:sivreal}
This formulation is obtained by separating real and imaginary parts of decision variables and constraints in the complex $(S,I,V)$-formulation from Sect.~\ref{s:sivform}. Decision variables are dealt with as follows:
\begin{eqnarray*}
  \forall b\in B \quad V_b &=& \re{V_b} +i\im{V_b} \\
  \forall (b,a,h)\in L\quad I_b &=& \re{I_{bah}} + i\im{I_{bah}} \\
  \forall (b,a,h)\in L\quad S_b &=& \re{S_{bah}} + i\im{S_{bah}} \\
  \forall g\in\mathscr{G} \quad \mathscr{S}_g &=& \re{\mathscr{S}_{g}} + i\im{\mathscr{S}_{g}}.
\end{eqnarray*}
As mentioned in Sect.~\ref{s:sivform}, the objective function is assumed to be real. We now tackle the constraints:
\begin{itemize}
\item Generated power bounds Eq.~\eqref{genpowerbounds}:
  \begin{eqnarray}
    \forall g\in\mathscr{G} \quad  \re{\underline{\mathscr{S}}_{g}}\le \re{\mathscr{S}_{g}} &\le& \re{\overline{\mathscr{S}}_{g}} \label{genpowerboundR} \\
    \forall g\in\mathscr{G} \quad  \im{\underline{\mathscr{S}}_{g}}\le \im{\mathscr{S}_{g}} &\le& \im{\overline{\mathscr{S}}_{g}}. \label{genpowerboundC}
  \end{eqnarray}
\item Bounds on the power magnitude Eq.~\eqref{powerbound}:
  \begin{equation}
    \forall (b,a,h)\in L\quad (\re{S_{bah}})^2 + (\im{S_{bah}})^2 \le \bar{S}_{bah}^2. \label{powerboundR}    
  \end{equation}
\item Bounds on phase differences Eq.~\eqref{phasediffbound2}-\eqref{phasediffboundaux}:
  \begin{eqnarray}
    \forall (b,a,h)\in L_0 && \tan(\underline{\eta}_{bah}) (\re{V_b}\re{V_a}+\im{V_b}\im{V_a}) \le \im{V_b}\re{V_a}-\re{V_b}\im{V_a} \label{phasediffbound1R} \\ \forall (b,a,h)\in L_0 && \im{V_b}\re{V_a}-\re{V_b}\im{V_a} \le \tan(\overline{\eta}_{bah}) (\re{V_b}\re{V_a}+\im{V_b}\im{V_a}) \label{phasediffbound2R}\\
    \forall (b,a,1)\in L_0 && \re{V_b}\re{V_a}+\im{V_b}\im{V_a} \ge 0. \label{phasediffboundauxR}
  \end{eqnarray}
\item Voltage bounds Eq.~\eqref{voltagebound}:
\begin{equation}
  \forall b\in B \quad \underline{V}_b^2 \le (\re{V_b})^2 + (\im{V_b})^2 \le \overline{V}_b^2. \label{voltageboundR}
\end{equation}
\item The reference bus constraints Eq.~\eqref{reference} are unchanged.
\item Power flow equations Eq.~\eqref{powerflow}:
  \begin{eqnarray}
    \forall b\in B\quad \sum\limits_{(b,a,h)\in L} \re{(S_{bah})} + \re{\tilde{S}_b} &=& -\re{A_b}|V_b|^2 + \sum\limits_{g\in\mathscr{G}_b} \re{\mathscr{S}_g} \label{powerflowR} \\
    \forall b\in B\quad \sum\limits_{(b,a,h)\in L} \im{(S_{bah})} + \im{\tilde{S}_b} &=& \im{A_b}|V_b|^2 + \sum\limits_{g\in\mathscr{G}_b} \im{\mathscr{S}_g}. \label{powerflowC}
  \end{eqnarray}
\item Power in terms of current Eq.~\eqref{powercurrent}:
  \begin{eqnarray}
    \forall (b,a,h)\in L \quad \re{(S_{bah})} &=& \re{V_b}\re{(I_{bah})} + \im{V_b}\im{(I_{bah})} \label{powercurrentR} \\
    \forall (b,a,h)\in L \quad \im{(S_{bah})} &=& \im{V_b}\re{(I_{bah})} - \re{V_b}\im{(I_{bah})}. \label{powercurrentC}
  \end{eqnarray}
\item Generalized Ohm's law Eq.~\eqref{ohmlaw1}-\eqref{ohmlaw2}:
  {\small
  \begin{eqnarray}
    \forall (b,a,h)\in L_0\quad \re{I_{bah}} &=& \re{(Y^{\mathsf{ff}}_{bah})}\re{V_b} - \im{(Y^{\mathsf{ff}}_{bah})}\im{V_b} + \re{(Y^{\mathsf{ft}}_{bah})}\re{V_a} - \im{(Y^{\mathsf{ft}}_{bah})}\im{V_a} \label{ohmlaw1R} \\
    \forall (b,a,h)\in L_0\quad \im{I_{bah}} &=& \re{(Y^{\mathsf{ff}}_{bah})}\im{V_b} + \im{(Y^{\mathsf{ff}}_{bah})}\re{V_b} + \re{(Y^{\mathsf{ft}}_{bah})}\im{V_a} + \im{(Y^{\mathsf{ft}}_{bah})}\re{V_a} \label{ohmlaw1C} \\
    \forall (b,a,h)\in L_0\quad \re{I_{abh}} &=& \re{(Y^{\mathsf{tf}}_{bah})}\re{V_b} - \im{(Y^{\mathsf{tf}}_{bah})}\im{V_b} + \re{(Y^{\mathsf{tt}}_{bah})}\re{V_a} - \im{(Y^{\mathsf{tt}}_{bah})}\im{V_a} \label{ohmlaw2R} \\
    \forall (b,a,h)\in L_0\quad \im{I_{abh}} &=& \re{(Y^{\mathsf{tf}}_{bah})}\im{V_b} + \im{(Y^{\mathsf{tf}}_{bah})}\re{V_b} + \re{(Y^{\mathsf{tt}}_{bah})}\im{V_a} + \im{(Y^{\mathsf{tt}}_{bah})}\re{V_a} \label{ohmlaw2C}.
  \end{eqnarray}
  }
\end{itemize}

\subsection{Cartesian voltage-only QCQP}
\label{s:vcart}
This formulation is derived from Eq.~\eqref{vonly2}; as such, it relies on magnitude bounds $\bar{I}$ on injected current rather than magnitude bounds $\bar{S}$ on injected power.

The objective function, linear in active power, is already a real function of real variables only. The bounds on generated power are enforced on real and complex parts separately, as in Eq.~\eqref{genpowerboundR}-\eqref{genpowerboundC}. We separate real and imaginary parts of the terms in the current magnitude bounds in Eq.~\eqref{Vcurrentbound1}-\eqref{Vcurrentbound2}, and obtain:
{\small
\begin{eqnarray}
  \forall (b,a,h)\in L_0 \quad |Y_{bah}^{\mathsf{ff}}|^2|V_b|^2 + |Y_{bah}^{\mathsf{ft}}|^2|V_a|^2 + 2\re{(Y_{bah}^{\mathsf{ff}} \conj{Y_{bah}^{\mathsf{ft}}} V_b \conj{V_a})} &\le& \bar{I}_{bah}^2 \label{currentbound1R} \\
  \forall (b,a,h)\in L_1 \quad |Y_{abh}^{\mathsf{tf}}|^2|V_a|^2 + |Y_{abh}^{\mathsf{tt}}|^2|V_b|^2 + 2\re{(Y_{abh}^{\mathsf{tf}} \conj{Y_{abh}^{\mathsf{tt}}} V_a \conj{V_b})}  &\le& \bar{I}_{bah}^2,\label{currentbound2R}
\end{eqnarray}
}%
where
{\small
  \begin{eqnarray}
  &&  |V_b|^2 = (\re{V_b})^2 + (\im{V_b})^2 \quad\land\quad
   |V_a|^2 = (\re{V_a})^2 + (\im{V_a})^2 \label{currentboundVV} \\ [0.3em]
  &&  \re{(Y_{bah}^{\mathsf{ff}} \conj{Y_{bah}^{\mathsf{ft}}} V_b \conj{V_a})} = \label{currentboundVba} \\
  &=& \re{Y_{bah}^{\mathsf{ff}}}\re{Y_{bah}^{\mathsf{ft}}} \re{V_b}\re{V_a} + \im{Y_{bah}^{\mathsf{ff}}}\im{Y_{bah}^{\mathsf{ft}}} \im{V_b}\im{V_a} + \re{Y_{bah}^{\mathsf{ff}}}\re{Y_{bah}^{\mathsf{ft}}} \im{V_b}\im{V_a} + \im{Y_{bah}^{\mathsf{ff}}}\im{Y_{bah}^{\mathsf{ft}}} \re{V_b}\re{V_a} \nonumber \\
  &-&   \re{Y_{bah}^{\mathsf{ff}}}\im{Y_{bah}^{\mathsf{ft}}} \re{V_b}\im{V_a} -
      \im{Y_{bah}^{\mathsf{ff}}}\re{Y_{bah}^{\mathsf{ft}}} \im{V_b}\re{V_a} + \re{Y_{bah}^{\mathsf{ff}}}\im{Y_{bah}^{\mathsf{ft}}} \im{V_b}\re{V_a} + \im{Y_{bah}^{\mathsf{ff}}}\re{Y_{bah}^{\mathsf{ft}}} \re{V_b}\im{V_a} \nonumber \\ [0.3em]
  &&  \re{(Y_{abh}^{\mathsf{tf}} \conj{Y_{abh}^{\mathsf{tt}}} V_a \conj{V_b})} =  \label{currentboundVab} \\
  &=& \re{Y_{abh}^{\mathsf{tf}}}\re{Y_{abh}^{\mathsf{tt}}} \re{V_a}\re{V_b} + \im{Y_{abh}^{\mathsf{tf}}}\im{Y_{abh}^{\mathsf{tt}}} \im{V_a}\im{V_b} + \re{Y_{abh}^{\mathsf{tf}}}\re{Y_{abh}^{\mathsf{tt}}} \im{V_a}\im{V_b} + \im{Y_{abh}^{\mathsf{tf}}}\im{Y_{abh}^{\mathsf{tt}}} \re{V_a}\re{V_b} \nonumber \\
  &-&   \re{Y_{abh}^{\mathsf{tf}}}\im{Y_{abh}^{\mathsf{tt}}} \re{V_a}\im{V_b} -
      \im{Y_{abh}^{\mathsf{tf}}}\re{Y_{abh}^{\mathsf{tt}}} \im{V_a}\re{V_b} + \re{Y_{abh}^{\mathsf{tf}}}\im{Y_{abh}^{\mathsf{tt}}} \im{V_a}\re{V_b} + \im{Y_{abh}^{\mathsf{tf}}}\re{Y_{abh}^{\mathsf{tt}}} \re{V_a}\im{V_b}.\nonumber
\end{eqnarray}
}%
The bounds on phase difference follow from Eq.~\eqref{phasediffbound2} using the identities:
\begin{eqnarray}
  \re{(V_b\conj{V_a})} = \re{V_b}\re{V_a} + \im{V_b}\im{V_a} &\qquad&
  \im{(V_b\conj{V_a})} = \im{V_b}\re{V_a} - \re{V_b}\im{V_a}. \label{vbva} 
\end{eqnarray}
We obtain:
\begin{eqnarray}
  \forall (b,a,h)\in L_0 && \tan(\underline{\eta}_{bah})(\re{V_b}\re{V_a} + \im{V_b}\im{V_a}) \le \im{V_b}\re{V_a} - \re{V_b}\im{V_a} \label{phasediffbound1VR} \\
  \forall (b,a,h)\in L_0 && \im{V_b}\re{V_a} - \re{V_b}\im{V_a} \le \tan(\overline{\eta}_{bah})(\re{V_b}\re{V_a} + \im{V_b}\im{V_a}) \label{phasediffbound2VR} \\
  \forall (b,a,1)\in L_0 && \re{V_b}\re{V_a} + \im{V_b}\im{V_a} \ge 0.\label{phasediffboundauxVR}
\end{eqnarray}
Voltage bounds are as in Eq.~\eqref{voltageboundR}, and reference bus constraints are as in Eq.~\eqref{reference}.

As concerns power, we separate real and imaginary parts of the decision variables in Eq.~\eqref{srepl1}-\eqref{srepl2}, and obtain:
{\small \begin{eqnarray}
  \forall (b,a,h) &\in& L_0: \nonumber \\
  \re{(S_{bah})} &=& \re{(Y_{bah}^{\mathsf{ff}})}|V_b|^2 + \re{(Y_{bah}^{\mathsf{ft}})}(\re{V_b}\re{V_a} + \im{V_b}\im{V_a}) + \im{(Y_{bah}^{\mathsf{ft}})}(\im{V_b}\re{V_a} - \re{V_b}\im{V_a}) \label{powercurr1VR} \\
  \im{(S_{bah})} &=& -\im{(Y_{bah}^{\mathsf{ff}})}|V_b|^2 + \re{(Y_{bah}^{\mathsf{ft}})}(\im{V_b}\re{V_a} - \re{V_b}\im{V_a}) - \im{(Y_{bah}^{\mathsf{ft}})}(\re{V_b}\re{V_a} + \im{V_b}\im{V_a}) \label{powercurr1VC} \\
  \forall (b,a,h) &\in& L_1: \nonumber \\
  \re{(S_{bah})} &=& \re{(Y_{abh}^{\mathsf{tt}})}|V_b|^2 + \re{(Y_{abh}^{\mathsf{tf}})}(\re{V_b}\re{V_a} + \im{V_b}\im{V_a}) + \im{(Y_{abh}^{\mathsf{tf}})}(\im{V_b}\re{V_a} - \re{V_b}\im{V_a}) \label{powercurr2VR} \\
  \im{(S_{bah})} &=&  - \im{(Y_{abh}^{\mathsf{tt}})}|V_b|^2 + \re{(Y_{abh}^{\mathsf{tf}})}(\im{V_b}\re{V_a} - \re{V_b}\im{V_a}) - \im{(Y_{abh}^{\mathsf{tf}})}(\re{V_b}\re{V_a} + \im{V_b}\im{V_a}), \label{powercurr2VC}
\end{eqnarray} }%
which can now be replaced in Eq.~\eqref{powerflowR}-\eqref{powerflowC} when written with separate sums over $L_0,L_1$ as in Eq.~\eqref{powerflow2}.

Note that we use the generalized Ohm's laws and the power definitions in terms of current only in order to operate replacement of $S_{bah}$ in the power flow equations Eq.~\eqref{powercurr1VR}-\eqref{powercurr2VC}.

\subsection{Polar formulation}
\label{s:polar}
The polar formulation of the ACOPF is obtained by the polar representation of complex voltage in terms of magnitude and phase:
\begin{equation}
  \forall b\in B \quad V_b = v_b e^{i\theta_b}, \label{Vpolar}
\end{equation}
where $v_b$ is the magnitude and $\theta_b$ is the phase (we remark that in Sect.~\ref{s:dyn2stat} we already introduced a scaled magnitude $V_b^{\max}=v_b/\sqrt{2}$). In the current setting, $v_b$ and $\theta_b$ are decision variables of the polar formulation. We also consider power generation variables $\mathscr{S}_g$ for $g\in\mathscr{G}$. 

Since we are describing a real (rather than complex) formulation, we write:
\begin{eqnarray}
  \forall b\in B \quad \re{V_b} &=& v_b \cos\theta_b  \label{VpolarR} \\
  \forall b\in B \quad \im{V_b} &=& v_b \sin\theta_b. \label{VpolarC}
\end{eqnarray}
Aside from Eq.~\eqref{VpolarR}-\eqref{VpolarC},  there is another implied relationship
\begin{equation}
  \forall b\in B \quad  v_b^2=|V_b|^2, \label{vV}
\end{equation}
between polar and cartesian formulations. All of these define nonconvex sets, however, so they are not exploited directly in MP formulations. By Eq.~\eqref{vV}, however, we can derive the bound constraints:
\begin{equation}
  \forall b\in B \quad  v_b\ge 0. \label{vnonneg}
\end{equation}
Moreover, by the periodicity of trigonometric functions, we can enforce the bound constraints:
\begin{equation}
  \forall b\in B \quad  -\pi\le \theta_b\le \pi. \label{phasebounds}
\end{equation}

The objective function is the same as in Eq.~\eqref{obj}. The same holds for the power generation bounds Eq.~\eqref{genpowerboundR}-\eqref{genpowerboundC}. The voltage magnitude bounds are
\begin{equation}
  \forall b\in B \quad \underline{V}_b \le v_b\le \overline{V}_b, \label{voltageboundvR}
\end{equation}
and the phase difference bounds are
\begin{equation}
  \forall (b,a,h)\in L_0 \quad \underline{\eta}_{bah} \le \theta_b-\theta_a\le\overline{\eta}_{bah}.\label{phasediffboundvR}
\end{equation}

We shall use the equations defining power in terms of current and voltage Eq.~\eqref{powercurrent} and the generalized Ohm's laws in order to write the injected power bounds Eq.~\eqref{powerbound} and the power flow equations Eq.~\eqref{powerflow}. We therefore have to express the injected power $S$ on the lines in function of the polar coordinate variables $v,\theta$. We can achieve this by replacing the right hand sides (rhs) of Eq.~\eqref{VpolarR}-\eqref{VpolarC} in the definitions of real and imaginary parts of voltage in Eq.~\eqref{powercurr1VR}-\eqref{powercurr2VC}, followed by the application of Ptolemy's identities:
\begin{eqnarray}
  \cos(\theta_b-\theta_a) &=& \cos\theta_b\cos\theta_a + \sin\theta_b\sin\theta_a \\
  \sin(\theta_b-\theta_a) &=& \sin\theta_b\cos\theta_a - \cos\theta_b\sin\theta_a.
\end{eqnarray}
This yields:
{\small \begin{eqnarray}
  \forall (b,a,h) &\in& L_0: \nonumber \\
  \re{(S_{bah})} &=& \re{(Y_{bah}^{\mathsf{ff}})}v_b^2  + \re{(Y_{bah}^{\mathsf{ft}})}v_bv_a\cos(\theta_b-\theta_a) + \im{(Y_{bah}^{\mathsf{ft}})}v_bv_a\sin(\theta_b-\theta_a) \label{powercurr1vR} \\
  \im{(S_{bah})} &=& -\im{(Y_{bah}^{\mathsf{ff}})}v_b^2  + \re{(Y_{bah}^{\mathsf{ft}})}v_bv_a\sin(\theta_b-\theta_a)  - \im{(Y_{bah}^{\mathsf{ft}})}v_bv_a\cos(\theta_b-\theta_a) \label{powercurr1vC} \\  
  \forall (b,a,h) &\in& L_1: \nonumber \\
  \re{(S_{bah})} &=& \re{(Y_{abh}^{\mathsf{tt}})}v_b^2 + \re{(Y_{abh}^{\mathsf{tf}})}v_bv_a\cos(\theta_b-\theta_a) + \im{(Y_{abh}^{\mathsf{tf}})}v_bv_a\sin(\theta_b-\theta_a) \label{powercurr2vR} \\
  \im{(S_{bah})} &=&  -\im{(Y_{abh}^{\mathsf{tt}})}v_b^2 + \re{(Y_{abh}^{\mathsf{tf}})}v_bv_a\sin(\theta_b-\theta_a) - \im{(Y_{abh}^{\mathsf{tf}})}v_bv_a\cos(\theta_b-\theta_a). \label{powercurr2vC}
\end{eqnarray} }%
As in Sect.~\ref{s:vcart}, these expressions can be used to replace injected power terms in Eq.~\eqref{powerflowR}-\eqref{powerflowC} when written with separate sums over $L_0,L_1$ as in Eq.~\eqref{powerflow2}. Unlike Sect.~\ref{s:vcart}, these expressions are also used to replace injected power terms in the power magnitude bounds Eq.~\eqref{powerboundR} (written as two separate constraints, quantified over $L_0$ and $L_1$). This yields two fourth-degree polynomial inequality constraints in $v$ which are linear in $\cos(\theta_b-\theta_a)$ and $\sin(\theta_b-\theta_a)$:
{\small\begin{eqnarray}
  \forall (b,a,h)\in L_0\ \left(\im{(Y_{bah}^{\mathsf{ff}})}\right)^{2} v_{b}^{4} + 2 \im{(Y_{bah}^{\mathsf{ff}})} \im{(Y_{bah}^{\mathsf{ft}})} v_{b}^{3}v_{a} \cos(\theta_b - \theta_a) && \nonumber \\
  - 2\im{(Y_{bah}^{\mathsf{ff}})} \re{(Y_{bah}^{\mathsf{ft}})} v_{b}^{3}v_{a} \sin(\theta_b - \theta_a) + \left(\im{(Y_{bah}^{\mathsf{ft}})}\right)^{2} v_{b}^{2}v_{a}^{2} && \nonumber \\
  + 2 \im{(Y_{bah}^{\mathsf{ft}})} \re{(Y_{bah}^{\mathsf{ff}})} v_{b}^{3}v_{a} \sin(\theta_b - \theta_a) + \left(\re{(Y_{bah}^{\mathsf{ff}})}\right)^{2} v_{b}^{4} && \nonumber \\
  + 2 \re{(Y_{bah}^{\mathsf{ff}})} \re{(Y_{bah}^{\mathsf{ft}})} v_{b}^{3}v_{a} \cos(\theta_b - \theta_a) + \left(\re{(Y_{bah}^{\mathsf{ft}})}\right)^{2} v_{b}^{2}v_{a}^{2} &\le& \bar{S}_{bah}^2 \label{powerbound1vR} \\
  \forall (b,a,h)\in L_1\ \left(\im{(Y_{abh}^{\mathsf{tf}})}\right)^{2} v_{b}^{2}v_{a}^{2} + 2 \im{(Y_{abh}^{\mathsf{tf}})} \im{(Y_{abh}^{\mathsf{tt}})} v_{b}^{3}v_{a} \cos(\theta_b - \theta_a) && \nonumber \\
  + 2 \im{(Y_{abh}^{\mathsf{tf}})} \re{(Y_{abh}^{\mathsf{tt}})} v_{b}^{3}v_{a} \sin(\theta_b - \theta_a) + \left(\im{(Y_{abh}^{\mathsf{tt}})}\right)^{2} v_{b}^{4} && \nonumber \\
  - 2 \im{(Y_{abh}^{\mathsf{tt}})} \re{(Y_{abh}^{\mathsf{tf}})} v_{b}^{3}v_{a} \sin(\theta_b - \theta_a) + \left(\re{(Y_{abh}^{\mathsf{tf}})}\right)^{2} v_{b}^{2}v_{a}^{2} && \nonumber \\ + 2 \re{(Y_{abh}^{\mathsf{tf}})} \re{(Y_{abh}^{\mathsf{tt}})} v_{b}^{3}v_{a}  \cos(\theta_b - \theta_a) + \left(\re{(Y_{abh}^{\mathsf{tt}})}\right)^{2} v_{b}^{4} &\le& \bar{S}_{bah}^2.\label{powerbound2vR}
\end{eqnarray}}%
We remark that the generated power variables $\mathscr{S}$ appear in Eq.~\eqref{powerflowR}-\eqref{powerflowC}, which, after the replacements mentioned above, are also part of this formulation.

\subsection{Jabr's relaxation}
\label{s:jabr}
A conic relaxation \cite{jabr} (also called Jabr's {\it radial relaxation}) can be obtained from the polar ACOPF formulation by replacing $v_bv_a\cos(\theta_b-\theta_a)$ by a new variable $c_{ba}$, and $v_bv_a\sin(\theta_b-\theta_a)$ by a new variable $s_{ba}$. More precisely, we define an index set $R=R_1\cup R_2$, where:
\begin{eqnarray*}
  R_1 &=& \{(b,b)\;|\;b\in B\} \\
  R_2 &=& \{(b,a)\;|\;(b,a,1)\in L\}.
\end{eqnarray*}
Now for all $(b,a)\in R$ we define new variables $c_{ba},s_{ba}$. Jabr's relaxation relies on $c,s,\mathscr{S}$ as decision variables.

If the following conditions
\begin{eqnarray}
  \forall (b,a)\in R \quad c_{ba} &=& v_bv_a\cos(\theta_b-\theta_a) \label{jabrexact1} \\
  \forall (b,a)\in R \quad s_{ba} &=& v_bv_a\sin(\theta_b-\theta_a) \label{jabrexact2}
\end{eqnarray}
held, then the Jabr's relaxation would turn out to be exact. From Eq.~\eqref{jabrexact1}-\eqref{jabrexact2} we infer:
\begin{eqnarray}
  \forall (b,a,1)\in L_0 \quad c_{ba} &=& c_{ab} \label{symm1J} \\
  \forall (b,a,1)\in L_0 \quad s_{ba} &=& -s_{ab} \label{symm2J} \\
  \forall b\in B \quad s_{bb} &=& 0 \label{szeroJ} \\
  \forall b\in B \quad c_{bb} &=& v_b^2 \label{jabr1} \\
  \forall (b,a,1)\in L_0 \quad c_{ba}^2 + s_{ba}^2 &=& v_b^2 v_a^2. \label{jabr2}
\end{eqnarray}
While the latter do not imply relaxation exactness, they are nonetheless valid constraints. Some of them are nonconvex, however, and two of them (Eq.~\eqref{jabr1}-\eqref{jabr2}) also involve the voltage magnitude variables $v$, which are not necessarily part of the relaxation. We therefore use Eq.~\eqref{jabr1} to replace $v$ in Eq.~\eqref{jabr2}, and relax the equality of Eq.~\eqref{jabr2} to a convex conic inequality:
\begin{equation}
  \forall (b,a,1)\in L_0 \quad c_{ba}^2 + s_{ba}^2 \le c_{bb}\, c_{aa}. \label{relaxJ}
\end{equation}
Moreover, by Eq.~\eqref{jabr1} we also have
\begin{equation}
  \forall b\in B\quad c_{bb}\ge 0. \label{cnonnegJ}
\end{equation}

We construct Jabr's relaxation as follows: the objective function is as in Eq.~\eqref{obj}, but we also assume it is convex quadratic. The power generation bounds Eq.~\eqref{genpowerboundR}-\eqref{genpowerboundC} are also part of the formulation. By Eq.~\eqref{jabr1}, the voltage magnitude bounds are
\begin{equation}
  \forall b\in B \quad \underline{V}_b \le c_{bb}\le \overline{V}_b. \label{voltageboundJ}
\end{equation}

By Eq.~\eqref{phasediffbound1VR}-\eqref{phasediffboundauxVR}, Eq.~\eqref{cnonnegJ}, and the fact that
\begin{eqnarray}
  \forall (b,a,1)\in L_0 \quad V_b\conj{V_a} &=& v_be^{i\theta_b}\,v_ae^{-i\theta_a} = v_bv_ae^{i(\theta_b-\theta_a)} \nonumber \\  &=& v_bv_a\cos(\theta_b-\theta_a)+iv_bv_a\sin(\theta_b-\theta_a) \nonumber \\
  &=& c_{ba} + i s_{ba}, \label{cartpol}
\end{eqnarray}
the phase difference bounds turn out to be:
\begin{eqnarray}
  \forall (b,a,h)\in L_0 && \tan(\underline{\eta}_{bah})\, c_{ba} \le s_{ba} \le \tan(\overline{\eta}_{bah})\, c_{ba} \label{phasediffboundJ1} \\
  \forall (b,a,h)\in L_0 && c_{ba} \ge 0. \label{phasediffboundJ2}
\end{eqnarray}

We obtain expression for the injected power by replacement of Eq.~\eqref{jabrexact1}-\eqref{jabrexact2} in Eq.~\eqref{powercurr1vR}-\eqref{powercurr2vC}:
\begin{eqnarray}
  \forall (b,a,h) &\in& L_0: \nonumber \\
  \re{(S_{bah})} &=& \re{(Y_{bah}^{\mathsf{ff}})}c_{bb} + \re{(Y_{bah}^{\mathsf{ft}})}c_{ba} + \im{(Y_{bah}^{\mathsf{ft}})}s_{ba} \label{powercurr1RJ} \\
  \im{(S_{bah})} &=& -\im{(Y_{bah}^{\mathsf{ff}})}c_{bb} + \re{(Y_{bah}^{\mathsf{ft}})}s_{ba} - \im{(Y_{bah}^{\mathsf{ft}})}c_{ba} \label{powercurr1CJ} \\
  \forall (b,a,h) &\in& L_1: \nonumber \\
  \re{(S_{bah})} &=& \re{(Y_{abh}^{\mathsf{tt}})}c_{bb} + \re{(Y_{abh}^{\mathsf{tf}})}c_{ba} + \im{(Y_{abh}^{\mathsf{tf}})}s_{ba} \label{powercurr2RJ} \\
  \im{(S_{bah})} &=&  - \im{(Y_{abh}^{\mathsf{tt}})}c_{bb} + \re{(Y_{abh}^{\mathsf{tf}})}s_{ba} - \im{(Y_{abh}^{\mathsf{tf}})}c_{ba}. \label{powercurr2CJ}
\end{eqnarray}
As in Sect.~\ref{s:vcart}, these expressions can be used to replace injected power terms in Eq.~\eqref{powerflowR}-\eqref{powerflowC} when written with separate sums over $L_0,L_1$ as in Eq.~\eqref{powerflow2}. Concerning the injected power bound inequalities Eq.~\eqref{powerboundR} (written as two separate constraints, quantified over $L_0$ and $L_1$), since Eq.~\eqref{powercurr1RJ}-\eqref{powercurr2CJ} are linear in $c,s$, we obtain the quadratic inequalities:
{\small\begin{eqnarray}
  \forall (b,a,h)\in L_0 \quad
  \left(\im{(Y_{bah}^{\mathsf{ff}})}\right)^2 c_{bb}^2 + 2\,\im{(Y_{bah}^{\mathsf{ff}})}\im{(Y_{bah}^{\mathsf{ft}})}c_{bb}c_{ba} && \nonumber \\
  - 2\,\im{(Y_{bah}^{\mathsf{ff}})}\re{(Y_{bah}^{\mathsf{ft}})}c_{bb}s_{ba} + \left(\re{(Y_{bah}^{\mathsf{ff}})}\right)^2 c_{bb}^2 && \nonumber \\
  + 2\,\re{(Y_{bah}^{\mathsf{ff}})}\im{(Y_{bah}^{\mathsf{ft}})} c_{bb}s_{ba} + 2\,\re{(Y_{bah}^{\mathsf{ff}})}\re{(Y_{bah}^{\mathsf{ft}})} c_{bb}c_{ba} && \nonumber \\
  + \left(\im{(Y_{bah}^{\mathsf{ft}})}\right)^2c_{bb}c_{aa} + \left(\re{(Y_{bah}^{\mathsf{ft}})}\right)^2c_{bb}c_{aa} &\le & \bar{S}_{bah}^2 \label{powerbound1J} \\
  \forall (b,a,h)\in L_1 \quad
  \left(\im{(Y_{abh}^{\mathsf{tf}})}\right)^2 c_{bb}c_{aa} + 2\,\im{(Y_{abh}^{\mathsf{tf}})}\im{(Y_{abh}^{\mathsf{tt}})}c_{bb}c_{ba} && \nonumber \\
  + 2\,\im{(Y_{abh}^{\mathsf{tf}})}\re{(Y_{abh}^{\mathsf{tt}})}c_{bb}s_{ba} + \left(\re{(Y_{abh}^{\mathsf{tf}})}\right)^2 c_{bb}c_{aa} && \nonumber \\
  - 2\,\re{(Y_{abh}^{\mathsf{tf}})}\im{(Y_{abh}^{\mathsf{tt}})} c_{bb}s_{ba} + 2\,\re{(Y_{abh}^{\mathsf{tf}})}\re{(Y_{abh}^{\mathsf{tt}})} c_{bb}c_{ba} && \nonumber \\
  + \left(\im{(Y_{abh}^{\mathsf{tt}})}\right)^2c_{bb}^2 + \left(\re{(Y_{abh}^{\mathsf{tt}})}\right)^2c_{bb}^2 &\le& \bar{S}_{bah}^2 \label{powerbound2J}.
\end{eqnarray}}%
We remark that Eq.~\eqref{symm1J}-\eqref{symm2J} are also part of the relaxation, and that the generated power variables $\mathscr{S}$ appear in Eq.~\eqref{powerflowR}-\eqref{powerflowC} with the replacements mentioned in Sect.~\ref{s:polar}.

\subsection{Mixed formulation}
\label{s:mixed}
An exact QCQP formulation can be derived from Jabr's conic relaxation by adding relationships between $c,s$ variables and rectangular $\re{V},\im{V}$ voltage variables. More precisely, the mixed formulation optimizes Eq.~\eqref{obj} subject to power generation bounds Eq.~\eqref{genpowerboundR}-\eqref{genpowerboundC}, voltage magnitude bounds Eq.~\eqref{cnonnegJ}-\eqref{voltageboundJ}, phase difference bounds Eq.~\eqref{phasediffboundJ1}-\eqref{phasediffboundJ2}, power flow equations obtained by replacement of Eq.~\eqref{powercurr1RJ}-\eqref{powercurr2CJ} into Eq.~\eqref{powerflowR}-\eqref{powerflowC} written as Eq.~\eqref{powerflow2}, power magnitude bounds Eq.~\eqref{powerbound1J}-\eqref{powerbound2J}, symmetry relations Eq.~\eqref{symm1J}-\eqref{symm2J}, as well as:
\begin{eqnarray}
\forall b\in B \quad c_{bb} &=& (\re{V_b})^2 + (\im{V_b})^2 \label{csVrel1} \\
\forall (b,a,1)\in L_0 \quad c_{ba} &=& \re{V_b}\re{V_a} + \im{V_b}\im{V_a} \label{csVrel2}\\
\forall (b,a,1)\in L_0 \quad s_{ba} &=& \im{V_b}\re{V_a} - \re{V_b}\im{V_a}. \label{csVrel3}
\end{eqnarray}

The exactness of the mixed formulation follows from the exactness of the polar formulation in Sect.~\ref{s:polar}, the relations Eq.~\eqref{jabrexact1}-\eqref{jabrexact2}, the fact Eq.~\eqref{cartpol}, and the identities Eq.~\eqref{vbva}.

\subsection{Matrix formulation}
\label{s:matrix}
In this section, we show a matrix formulation of the voltage-only formulation presented above. This formulation is the one usually developed in solvers and it is inspired on the ones found in \cite{lavaei,molzahnsurvey}.

Rearranging terms, from Eq.~\eqref{vonly} we have that for each bus $b\in B$,
\begin{eqnarray}
\sum\limits_{g\in\mathscr{G}_b} \mathscr{S}_g - \tilde{S}_b
 &=&
\conj{A_b}|V_b|^2 + \sum\limits_{(b,a,h)\in L_0} (\conj{Y_{bah}^{\mathsf{ff}}} |V_b|^2 + \conj{Y_{bah}^{\mathsf{ft}}} V_b\conj{V_a})\ + \notag \\
&& + \sum\limits_{(b,a,h)\in L_1}(\conj{Y_{abh}^{\mathsf{tf}}} V_b\conj{V_a} + \conj{Y_{abh}^{\mathsf{tt}}} |V_b|^2).  \label{powbal}
\end{eqnarray}
The left hand side (lhs) of Eq.~\eqref{powbal} equals the net complex power injected to bus $b$. This equation can be rewritten as follows:
\begin{eqnarray}
\sum\limits_{g\in\mathscr{G}_b} \mathscr{S}_g - \tilde{S}_b
 &=& V_b\bigg[
\conj{A_b}\conj{V_b} + \sum\limits_{(b,a,h)\in L_0} (\conj{Y_{bah}^{\mathsf{ff}}} \conj{V_b} + \conj{Y_{bah}^{\mathsf{ft}}} \conj{V_a})\ + \notag \\
&& \quad + \sum\limits_{(b,a,h)\in L_1}(\conj{Y_{abh}^{\mathsf{tf}}} \conj{V_a} + \conj{Y_{abh}^{\mathsf{tt}}} \conj{V_b})\bigg].  \label{powbal2}
\end{eqnarray}

In section \ref{s:linemodel} we defined the admittance matrix $\mathbf{Y}_{ba}$ of the line $(b,a)$. Now, we define the admitance matrix of the network, let $\mathcal{Y}\in\mathbb{C}^{n\times n}$ be defined by
\begin{eqnarray*}
  \forall b \in B\quad\mathcal{Y}_{bb}
  &=& A_b + \sum\limits_{(b,a,h)\in L_0} Y_{bah}^{\mathsf{ff}} + \sum\limits_{(b,a,h)\in L_1} Y_{abh}^{\mathsf{tt}}\\
  \forall b,a \in B,b\neq a \quad \mathcal{Y}_{ba}
  &=& \sum\limits_{(b,a,h)\in L_0} Y_{bah}^{\mathsf{ft}} + \sum\limits_{(b,a,h)\in L_1} Y_{abh}^{\mathsf{tf}}.
\end{eqnarray*}
Denote by $e_1,\ldots,e_n$ the standard basis vectors in $\mathbb{R}^n$.
Let $\Psi_b:=e_b\transpose{e_b}\mathcal{Y}$, define the voltage vector $V := \transpose{(V_1\ \cdots\ V_n)}$ and
the matrix $W:=
{\small\left[\begin{array}{c} \re{V}\\ \im{V}\end{array}\right]
\transpose{\left[\begin{array}{c} \re{V}\\ \im{V}\end{array}\right]}}$.
Therefore, Eq.~\eqref{powbal2} can be expressed as:
\begin{eqnarray}
   \notag\sum\limits_{g\in\mathscr{G}_b} \mathscr{S}_g - \tilde{S}_b
   &=& V_b\sum_{a\in B} \conj{\mathcal{Y}_{ba}}\conj{V_a} = V_b \conj{(\transpose{e_b}\mathcal{Y}V)}\\
   \notag &=& (\transpose{e_b}V)(\hermitian{V}\hermitian{\mathcal{Y}}e_b)
   = (\hermitian{V}\hermitian{\mathcal{Y}}e_b)(\transpose{e_b}V)\\
   &=& \conj{(\hermitian{V}\Psi_bV)} \label{lavpf}\\
   \notag&=& \conj{\big(\transpose{(\re{V}-i\im{V})}(\re{\Psi}_b+i\im{\Psi}_b)(\re{V}+i\im{V})\big)}\\
   \notag&=& 
   (\transpose{\re{V}}\re{\Psi}_b\re{V}+\transpose{\im{V}}\re{\Psi}_b\im{V}
   +\transpose{\im{V}}\im{\Psi}_b\re{V}-\transpose{\re{V}}\im{\Psi}_b\im{V})\\
   \notag&&+\ i
   (\transpose{\im{V}}\re{\Psi}_b\re{V}-\transpose{\re{V}}\re{\Psi}_b\im{V}
   -\transpose{\re{V}}\im{\Psi}_b\re{V}-\transpose{\im{V}}\im{\Psi}_b\im{V})\\
   \notag&=&
   \transpose{\left[\begin{array}{c} \re{V}\\ \im{V}\end{array}\right]}
   \left[\begin{array}{ccc}\re{\Psi}_b && -\im{\Psi}_b \\
         \im{\Psi}_b && \re{\Psi}_b\end{array}\right]
   \left[\begin{array}{c} \re{V}\\ \im{V}\end{array}\right]
   - i
   \transpose{\left[\begin{array}{c} \re{V}\\ \im{V}\end{array}\right]}
   \left[\begin{array}{ccc}\im{\Psi}_b && \re{\Psi}_b \\
         -\re{\Psi}_b && \im{\Psi}_b\end{array}\right]
   \left[\begin{array}{c} \re{V}\\ \im{V}\end{array}\right]\\
   \notag&=&
   \frac{1}{2}\transpose{\left[\begin{array}{c} \re{V}\\ \im{V}\end{array}\right]}
   \left[\begin{array}{ccc}\re{\Psi}_b+\transpose{\re{\Psi}_b} && \transpose{\im{\Psi}_b}-\im{\Psi}_b \\
         \im{\Psi}_b-\transpose{\im{\Psi}_b} && \re{\Psi}_b+\transpose{\re{\Psi}_b}\end{array}\right]
   \left[\begin{array}{c} \re{V}\\ \im{V}\end{array}\right]\\
   \notag&& -\ i\frac{1}{2}
   \transpose{\left[\begin{array}{c} \re{V}\\ \im{V}\end{array}\right]}
   \left[\begin{array}{ccc}\im{\Psi}_b+\transpose{\im{\Psi}_b} && \re{\Psi}_b-\transpose{\re{\Psi}_b} \\
         \transpose{\re{\Psi}_b}-\re{\Psi}_b && \im{\Psi}_b+\transpose{\im{\Psi}_b}\end{array}\right]
   \left[\begin{array}{c} \re{V}\\ \im{V}\end{array}\right].
\end{eqnarray}
In consequence, the power flow balance equations are:
\begin{eqnarray*}
 \sum\limits_{g\in\mathscr{G}_b} \re{\mathscr{S}_g} - \re{\tilde{S}_b}
 &=& \trace{\left(\frac{1}{2}
   \left[\begin{array}{ccc}\re{\Psi}_b+\transpose{\re{\Psi}_b} && \transpose{\im{\Psi}_b}-\im{\Psi}_b \\
         \im{\Psi}_b-\transpose{\im{\Psi}_b} && \re{\Psi}_b+\transpose{\re{\Psi}_b}\end{array}\right]
   \left[\begin{array}{c} \re{V}\\ \im{V}\end{array}\right]
   \transpose{\left[\begin{array}{c} \re{V}\\ \im{V}\end{array}\right]}\right)}
 \\
 \sum\limits_{g\in\mathscr{G}_b} \im{\mathscr{S}_g} - \im{\tilde{S}_b}
 &=& \trace{\left(\frac{-1}{2}
   \left[\begin{array}{ccc}\im{\Psi}_b+\transpose{\im{\Psi}_b} && \re{\Psi}_b-\transpose{\re{\Psi}_b} \\
         \transpose{\re{\Psi}_b}-\re{\Psi}_b && \im{\Psi}_b+\transpose{\im{\Psi}_b}\end{array}\right]
   \left[\begin{array}{c} \re{V}\\ \im{V}\end{array}\right]
   \transpose{\left[\begin{array}{c} \re{V}\\ \im{V}\end{array}\right]}\right)}.
\end{eqnarray*}
Finally, we write
\begin{eqnarray} \label{trace-bal}
  \sum\limits_{g\in\mathscr{G}_b} \mathscr{S}_g &=& \tilde{S}_b +
  \trace{(\mathbf{\Psi}_b W)} + i \trace{(\mathbf{\hat{\Psi}}_b  W)}
\end{eqnarray}
where 
\begin{equation*}
\mathbf{\Psi}_b := \dfrac{1}{2}
   \left[\begin{array}{ccc}\re{\Psi}_b+\transpose{\re{\Psi}_b} && \transpose{\im{\Psi}_b}-\im{\Psi}_b \\
         \im{\Psi}_b-\transpose{\im{\Psi}_b} && \re{\Psi}_b+\transpose{\re{\Psi}_b}\end{array}\right],
\qquad
\mathbf{\hat{\Psi}}_b := \dfrac{-1}{2}
   \left[\begin{array}{ccc}\im{\Psi}_b+\transpose{\im{\Psi}_b} && \re{\Psi}_b-\transpose{\re{\Psi}_b} \\
         \transpose{\re{\Psi}_b}-\re{\Psi}_b && \im{\Psi}_b+\transpose{\im{\Psi}_b}\end{array}\right].
\end{equation*}

On the other hand, the flow on a branch $(b,a,h)\in L$ can be written as follows. Define
\begin{equation*}
  \Phi_{bah} := \begin{cases}
     Y_{bah}^{\mathsf{ff}} e_b\transpose{e_b} + Y_{bah}^{\mathsf{ft}} e_b\transpose{e_a},&\mbox{if }(b,a,h)\in L_0\\
     Y_{abh}^{\mathsf{tf}} e_b\transpose{e_a} + Y_{abh}^{\mathsf{tt}} e_b\transpose{e_b},&\mbox{if }(b,a,h)\in L_1.
  \end{cases}
\end{equation*}
Then,
\begin{equation*}
  S_{bah} = \hermitian{V} \hermitian{\Phi_{bah}}V
  = \conj{(\hermitian{V} \Phi_{bah}V)},
\end{equation*}
and, by following the steps in Eq.~\eqref{lavpf}, we have
\begin{eqnarray}  \label{trace-pow}
 S_{bah} &=&
  \trace{(\mathbf{\Phi}_{bah} W)} + i\trace{(\mathbf{\hat{\Phi}}_{bah} W)}\\
 \notag \mbox{where }\qquad \mathbf{\Phi}_{bah} &:=& \dfrac{1}{2}
   \left[\begin{array}{ccc}\re{\Phi}_{bah}+\transpose{\re{\Phi}_{bah}} && \transpose{\im{\Phi}_{bah}}-\im{\Phi}_{bah} \\
         \im{\Phi}_{bah}-\transpose{\im{\Phi}_{bah}} && \re{\Phi}_{bah}+\transpose{\re{\Phi}_{bah}}\end{array}\right]\\
 \notag \mbox{and }\qquad  \mathbf{\hat{\Phi}}_{bah} &:=& \dfrac{-1}{2}
   \left[\begin{array}{ccc}\im{\Phi}_{bah}+\transpose{\im{\Phi}_{bah}} && \re{\Phi}_{bah}-\transpose{\re{\Phi}_{bah}} \\
         \transpose{\re{\Phi}_{bah}}-\re{\Phi}_{bah} && \im{\Phi}_{bah}+\transpose{\im{\Phi}_{bah}}\end{array}\right].
\end{eqnarray}
Therefore, thanks to Shur's complement formula, the power flow constraint $|S_{bah}|^2=\trace{(\mathbf{\Phi}_{bah} W)}^2 + \trace{(\mathbf{\hat{\Phi}}_{bah} W)}^2\leq \bar{S}_{bah}^2$ can be written as the PSD constraint
\begin{equation} \label{trace-pow-sdp}
 \left[\begin{array}{ccccc} -\bar{S}_{bah}^2 && \trace{(\mathbf{\Phi}_{bah}W)} && \trace{(\mathbf{\hat{\Phi}}_{bah}W)}\\
 \trace{(\mathbf{\Phi}_{bah}W)} && -1 && 0\\
 \trace{(\mathbf{\hat{\Phi}}_{bah}W)} && 0 && -1
 \end{array}\right]\preceq 0.
\end{equation}

Moreover, since $V_b\conj{V_a}=\conj{(\hermitian{V}e_b\transpose{e_a}V)}$,
\begin{eqnarray} \label{trace-VV}
 V_b\conj{V_a} &=& \trace{(\mathbf{\Theta}_{ba} W)} + i\trace{(\mathbf{\hat{\Theta}}_{ba} W)}\\
\notag  \mbox{where }\qquad \mathbf{\Theta}_{ba} &:=& \dfrac{1}{2}
   \left[\begin{array}{ccc}e_b\transpose{e_a}+e_a\transpose{e_b} && 0 \\
         0 && e_b\transpose{e_a}+e_a\transpose{e_b}\end{array}\right]\\
\notag \mbox{and }\qquad  \mathbf{\hat{\Theta}}_{ba} &:=& \dfrac{-1}{2}
   \left[\begin{array}{ccc}0 && e_b\transpose{e_a}-e_a\transpose{e_b} \\
         e_a\transpose{e_b}-e_b\transpose{e_a} && 0\end{array}\right].
\end{eqnarray}

Using the representation of the balance equations, the power flow and the multiplication of voltages expressed in Eqs.~\eqref{trace-bal}, \eqref{trace-pow} and \eqref{trace-VV}, we obtain the following formulation of the ACOPF problem presented in the previous sections, over the variables $\re{V},\im{V}\in\mathbb{R}^n$, $W\in\mathbb{R}^{2n\times 2n}$ and $\mathscr{S}_g\in\mathbb{C}$ for each $g\in\mathscr{G}$:
\begin{equation}
  \left.
  \begin{array}{rrcl} 
    \min & \sum\limits_{g\in\mathscr{G}} (c_{g1}\re{\mathscr{S}_g} + c_{g0}) && \\
    \forall g\in \mathscr{G} &  \underline{\mathscr{S}}_{g}\
      \leq\ \mathscr{S}_{g}\ &\leq&\ \overline{\mathscr{S}}_{g}  \\ [0.4em]
    \forall b\in B &  \tilde{S}_b + \trace{(\mathbf{\Psi}_bW)} + i\trace{(\mathbf{\hat{\Psi}}_bW)}
        &=& \sum\limits_{g\in\mathscr{G}_b} {\mathscr{S}_g}\\ 
    \forall (b,a,h)\in L &
        \trace{(\mathbf{\Phi}_{bah} W)}^2 + \trace{(\mathbf{\hat{\Phi}}_{bah} W)}^2
        &\leq& \bar{S}^2_{bah} \\ [0.4em]
    \forall (b,a,h)\in L_0 & \quad
      [\tan(\underline{\eta}_{bah}),\tan(\overline{\eta}_{bah})]\trace{(\mathbf{\Theta}_{ba} W)}
      &\ni & \trace{(\mathbf{\hat{\Theta}}_{ba} W)}  \\ [0.4em]
    \forall (b,a,h)\in L_0 & \trace{(\mathbf{\Theta}_{ba} W)}
      &\geq & 0  \\ [0.4em]
    & \transpose{e_r}\im{V} = 0 \quad \land\quad \transpose{e_r}\re{V} &\ge & 0 \\ [0.4em]
    \forall b\in B & \underline{V}_b^2\ \le\ \trace{(\mathbf{\Theta}_{bb} W)}  &\le& \overline{V}_b^2 \\ [0.1em]
    & {\small\left[\begin{array}{c} \re{V}\\ \im{V}\end{array}\right]
      \transpose{\left[\begin{array}{c} \re{V}\\ \im{V}\end{array}\right]}}
      & = & W.
  \end{array}
  \right\} \label{lavaei-acopf}
\end{equation}
The reference bus constraint $\transpose{e_r}\im{V}=0$ can be expressed in terms of $W$ as
${\small\transpose{\left[\begin{array}{c} 0\\ e_r\end{array}\right]}}W
{\small\left[\begin{array}{c} 0\\ e_r\end{array}\right]}=0$.
Note that if $(\re{V},\im{V},W)$ is a feasible point of Eq.~\eqref{lavaei-acopf}, then $(-\re{V},-\im{V},W)$ is also feasible (with the exception of constraint $\transpose{e_r}\re{V}\geq0$) with the same objective function value. Therefore, we can omit constraint $\transpose{e_r}\re{V}\geq0$ and choose the corresponding solution that satisfy this inequality.

We obtain a real SDP relaxation of Eq.~\eqref{lavaei-acopf} by relaxing the constraint $W=
{\small\left[\begin{array}{c} \re{V}\\ \im{V}\end{array}\right]
\transpose{\left[\begin{array}{c} \re{V}\\ \im{V}\end{array}\right]}}$ to
$W\succeq_{\mathsf{H}} 0$ and, therefore, we can omit the voltage variables in the formulation:
\begin{equation}
  \left.
  \begin{array}{rrcl} 
    \min & \sum\limits_{g\in\mathscr{G}} (c_{g1}\re{\mathscr{S}_g} + c_{g0}) && \\
    \forall g\in \mathscr{G} &  \underline{\mathscr{S}}_{g}\
      \leq\ \mathscr{S}_{g}\ &\leq&\  \overline{\mathscr{S}}_{g}  \\ [0.4em]
    \forall b\in B &  \tilde{S}_b + \trace{(\mathbf{\Psi}_bW)} + i\trace{(\mathbf{\hat{\Psi}}_bW)}
        &=& \sum\limits_{g\in\mathscr{G}_b} {\mathscr{S}_g}\\ 
    \forall (b,a,h)\in L & {\small
       \left[\begin{array}{ccc} -\bar{S}_{bah}^2 & \trace{(\mathbf{\Phi}_{bah}W)} & \trace{(\mathbf{\hat{\Phi}}_{bah}W)}\\
       \trace{(\mathbf{\Phi}_{bah}W)} & -1 & 0\\
       \trace{(\mathbf{\hat{\Phi}}_{bah}W)} & 0 & -1
       \end{array}\right]} &\preceq& 0 \\ [1.2em]
    \forall (b,a,h)\in L_0 & \quad
      [\tan(\underline{\eta}_{bah}),\tan(\overline{\eta}_{bah})]\trace{(\mathbf{\Theta}_{ba} W)}
      &\ni & \trace{(\mathbf{\hat{\Theta}}_{ba} W)}  \\ [0.4em]
    \forall (b,a,h)\in L_0 & \trace{(\mathbf{\Theta}_{ba} W)}
      &\geq & 0  \\ [0.4em]
    & \trace{\big({\small\left[\begin{array}{c} 0\\ e_r\end{array}\right]\transpose{\left[\begin{array}{c} 0\\ e_r\end{array}\right]}} W\big)} &= & 0 \\ [0.4em]
    \forall b\in B & \underline{V}_b^2\ \le\ \trace{(\mathbf{\Theta}_{bb} W)} &\le& \overline{V}_b^2 \\ [0.4em]
    & W & \succeq_{\mathsf{H}} & 0.
  \end{array}
  \right\} \label{lavaei-sdp}
\end{equation}

\section{Recent results on ACOPF relaxations}
\label{s:rel}

\LEO{In this section we review a few important papers in the literature concerning relaxations of the ACOPF, notably using SDP and derived techniques.}

\subsection{Zero duality gaps in SDP relaxations}
\label{s:lavaeilow}
The authors of \cite{lavaei} consider the dual of a real SDP relaxation of the ACOPF. Under particular conditions, the duality gap between primal and dual formulations is shown to be zero. A solution of the original OPF problem can therefore be recovered from the dual optimal solution. The dual SDP considered in \cite{lavaei} is similar to Eq.~\eqref{lavaei-sdp}. Table \ref{tab:lavaei-notation} shows the notational differences between \cite{lavaei} and Eq.~\eqref{lavaei-sdp}. There are four differences between these two formulations: (i) parallel lines are not considered in \cite{lavaei}; (ii) at most one generator is assumed to be attached to each bus in \cite{lavaei}, which allows the power balance equations to be replaced in the generated power bound inequalities; (iii) the power magnitude bound on lines in \cite{lavaei} is only imposed on the real part of the power variables, i.e.~$|\re{S}_{bah}|\leq \overline{P}_{bah}$; (iv) the phase difference constraints are replaced by $|V_b-V_a|\leq \Delta_{ba}$, where $\Delta_{ba}$ is a given parameter.
\begin{table}
\begin{center}
\caption{Equivalence of notations between Eq.~\eqref{lavaei-sdp}, \cite{lavaei} and \cite{molzahn2013}.}
\label{tab:lavaei-notation}
\begin{tabular}{c|c|c|p{3.5cm}}
Our & Notation & Notation & \\
notation & in \cite{lavaei}  & in \cite{molzahn2013} &Comments\\ \hline
$b$ & $k$ & $i$ & generic bus\\
$(b,a,h)$ & $(l,m)$ & $k=(k_l,k_m)$ & generic line \\ \hline
$V$ & $\mathbf{V}$ & & vector of complex voltages\\
$W$ & $W$ & $\mathbf{W}$ & \\
$\sum\limits_{g\in\mathscr{G}_b} \mathscr{S}_g$ &  $P_{G_k}+iQ_{G_k}$ & $\sum\limits_{g\in{\cal G}_i}(P_{G_g}+iQ_{G_g})$ &  generation at a generic bus\\
$\tilde{S}_b$ & $P_{D_k}+iQ_{D_k}$ & $P_{D_i}+iQ_{D_i}$ &  load at a generic bus\\
$\mathcal{Y}$ & $Y$ & $\mathbf{Y}$ & admittance matrix\\
$\Psi_b$ & $Y_k$ & $Y_i$ & $\Psi_b=e_b\transpose{e_b}\mathcal{Y}$\\ \hline \hline
\multicolumn{3}{c|}{Constraint matrices}& Related constraint\\ \hline
$\mathbf{\Psi}_b$ & $\mathbf{Y}_k$  & $\mathbf{Y}_i$ & real power injection\\
$\mathbf{\hat{\Psi}}_b$ & $\mathbf{\bar{Y}}_k$ & $\mathbf{\bar{Y}}_i$  & complex power injection\\
$\mathbf{\Theta}_{bb}$ & $M_k$ & $\mathbf{M}_i$ & voltage magnitude\\
$\mathbf{\Phi}_{bah}$ & $\mathbf{Y}_{lm}$  &
{$\small \begin{cases}\mathbf{Z}_{k_l},&(b,a,h)\in L_0\\ \mathbf{Z}_{k_m},&(b,a,h)\in L_1\end{cases}$} & branch active power flow\\
$\mathbf{\hat{\Phi}}_{bah}$ & $\mathbf{\bar{Y}}_{lm}$ &
{$\small \begin{cases}\mathbf{\bar{Z}}_{k_l},&(b,a,h)\in L_0\\ \mathbf{\bar{Z}}_{k_m},&(b,a,h)\in L_1\end{cases}$}  & branch reactive power flow
\end{tabular}
\end{center}
\end{table}

The algorithm proposed in \cite{lavaei} to obtain an optimal solution of the ACOPF is based on the zero duality gap between the ACOPF and the dual SDP relaxation: this occurs if the matrix
\begin{align*}
  A(x,r):=&\sum_{b\in B}(x_{b,1}\mathbf{\Psi}_b+x_{b,2}\mathbf{\hat{\Psi}}_b+x_{b,3}\mathbf{\Theta}_{bb})\ +\\
  &+\sum_{(b,a,h)\in L} \big((2r_{bah,1}+x_{bah,1})\mathbf{\Phi}_{bah} + 2r_{bah,2}\mathbf{\hat{\Phi}}_{bah}+\\
  &\hspace{1.8cm}+x_{bah,2}(\mathbf{\Theta}_{bb}+\mathbf{\Theta}_{aa}-2\mathbf{\Theta}_{ba})\big)
\end{align*}
found in the dual constraint $A(x,r)\succeq_{\mathsf{H}} 0$ (where $x\geq0$ and $r$ are the Lagrange multipliers associated to the primal constraints) has a zero eigenvalue of multiplicity at most 2, when the matrix is evaluated in the dual optimal solution $(x^{\text{opt}},r^{\text{opt}})$. 
Therefore, a globally optimal solution of the ACOPF can be recovered in polynomial time from a nonzero vector in the null space of the aforementioned matrix. 

Similar conditions for ensuring zero duality gap for DC power distribution networks are also given in \cite{lavaei}. Moreover, it is shown that such conditions happen almost always in DC networks. Empirical results are shown to prove the efficiency of this method on IEEE benchmark systems with 14, 30, 57, 118, and 300 buses.

\subsection{Applicability to real-life cases}
\label{s:molzahn}
The applicability of the ideas given in \cite{lavaei} to real-life cases was carried out in a sequence of papers co-authored by Molzahn {\it et al.}. 

The generalizations put forth in \cite{molzahn2014-3} are as follows: (i) non-zero line resistances are allows; (ii) conditions for zero duality gap are derived from the KKT conditions: given a solution $V$ of the dual SDP relaxation in terms of voltage values, the corresponding rank one matrix $W$, and a dual feasible matrix $A\succeq_{\mathsf{H}} 0$, the value of the dual SDP relaxation is globally optimal for the ACOPF if $A,W$ satisfy the slack complementarity condition $\trace{(AW)}=0$.

In \cite{molzahn2013}, the SDP formulations of \cite{lavaei} are extended to the case of more than one generator connected to each bus, and multiple (parallel) lines between buses. Moreover, the generation cost may be a quadratic or piece-wise linear function of the real power generated. However, phase difference and reference voltage fixing constraints are not considered. The dual SDP relaxation presented in \cite{molzahn2013} integrates these generalizations. An extension of the PSD matrix completion theorem in \cite{jabr2012} allows the replacement of the PSD variable matrix $W$ (as mentioned above) by multiple but smaller-sized PSD variable matrices, at the cost of some additional linking constraints between the components of the smaller matrices. The size reduction extent depends on the sparsity of the power network $G$. The solution algorithm proposed in \cite{molzahn2013} takes care of the trade-off between size reduction and the additional linking constraints. Numerical results show the efficiency of this method in large networks, such as the IEEE 300-bus system and the 3012-bus model of the Polish system.

We note that Table \ref{tab:lavaei-notation} also shows the notational differences with \cite{molzahn2013}. To help with readability, we also note that \cite[Prob.~(2)]{molzahn2013} and \cite[Prob.~(10)]{molzahn2014-3} use the symbols $S_k,S_{lm}$ to mean injected power on a line. 

\subsection{Use of Lasserre's relaxation hierarchies}
As mentioned in Sect.~\ref{s:sdprel}, and in particular in Eq.~\eqref{vonly2}, the ACOPF can be formulated as a QCQP, which is a Polynomial Programming (PP) problem of degree 2. Such formulations were investigated in \cite{ghaddar2016,zuluaga-acopf}. We report the formulation of interest using the notation introduced in Sect.~\ref{s:matrix}, with variables $\re{\mathscr{S}_g}$, $\re{S_{bah}}$, $\im{S_{bah}}$ and $x=\big[\transpose{\re{V}}\ \transpose{\im{V}}\big]$:
\begin{equation}
  \left.
  \begin{array}{rrcl} 
    \min & \sum\limits_{g\in\mathscr{G}} (c_{g2}(\re{\mathscr{S}_g})^2+ c_{g1}(\trace{(\mathbf{\Psi}_gx\transpose{x})} + \re{\tilde{S}_g}) &+& c_{g0}) \\
    \forall b\in B &  \re{\underline{\mathscr{S}}_{b}}\
      \leq\ \trace{(\mathbf{\Psi}_bx\transpose{x})} + \re{\tilde{S}_b}
      &\leq& \re{\overline{\mathscr{S}}_{b}}  \\ [0.4em]
    \forall b\in B &  \im{\underline{\mathscr{S}}_{b}}\
      \leq\ \trace{(\mathbf{\hat{\Psi}}_bx\transpose{x})} + \im{\tilde{S}_b}
      &\leq& \im{\overline{\mathscr{S}}_{b}}  \\ [0.4em]
    \forall b\in B & \underline{V}_b^2\ \le\ \trace{(\mathbf{\Theta}_{bb} x\transpose{x})} 
      &\le& \overline{V}_b^2 \\ [0.4em]
    \forall (b,a,h)\in L &
        (\re{S_{bah}})^2 + (\im{S_{bah}})^2  &\leq& \bar{S}^2_{bah} \\ [0.4em]
    \forall g \in \mathscr{G}& \trace{(\mathbf{\Psi}_g x\transpose{x})} + \re{\tilde{S}_g}
        &=& \re{\mathscr{S}_g} \\ [0.4em]
    \forall (b,a,h)\in L &  \trace{(\mathbf{\Phi}_{bah} x\transpose{x})} &=& \re{S_{bah}} \\ [0.4em]
    \forall (b,a,h)\in L &  \trace{(\mathbf{\hat{\Phi}}_{bah} x\transpose{x})} &=& \im{S_{bah}}.
  \end{array}
  \right\} \label{pp2-acopf}
\end{equation}
We note that \cite{ghaddar2016,zuluaga-acopf} assume that there are no parallel arcs, and that there is at most one generator attached to each bus. Moreover, the power flow equations are implicit in the power generation bounds, through the replacement $\underline{\mathscr{S}}_{b}=\sum_{g\in\mathscr{G}_b}\underline{\mathscr{S}}_{g}$ and $\overline{\mathscr{S}}_{b}=\sum_{g\in\mathscr{G}_b}\overline{\mathscr{S}}_{g}$, where the sums are equal to zero if $\mathscr{G}_b=\varnothing$.

The solution approaches in \cite{ghaddar2016,zuluaga-acopf} are both based on the well-known Lasserre relaxation hierarchies \cite{lasserre}, which builds a sequence of SDP relaxations that approximate the dual of a PP such as Eq.~\eqref{pp2-acopf}. The decision variables of these SDP relaxations are the coefficients of the monomials appearing in Sum-Of-Squares (SOS) polynomials up to certain degree. A sequence of SDP relaxations is obtained by increasing the degree of the SOS polynomials, which in turn increases the tightness of the approximation.

The first level of the Lasserre hierarchy that approximates Eq.~\eqref{pp2-acopf} is shown by \cite{ghaddar2016} to be equivalent to the dual SDP relaxation in \cite{lavaei} mentioned above. Linear and SOCP approximations of the SDP cone based on Diagonal Dominance (DD) were employed in \cite{zuluaga-acopf} for more efficient computation, as explained in \cite{ahmadi2014,ahmadimajumdar}.

The sparse structure of the polynomials involved in Eq.~\eqref{pp2-acopf} is exploited in \cite{zuluaga-acopf} in order to reduce the number of monomials appearing in the DD approximations of Lasserre's hierarchy (each approximating formulation in this hierarchy is called a ``structured PP-SDSOS${}_r$''). We note that this structure is also studied in \cite{bose2015}. The first level of the SOCP DD approximation of Lasserre's hierarchy is shown in \cite{zuluaga-acopf} to be equivalent to the dual of the SOCP relaxation of the ACOPF presented in \cite{bose2015}:
\begin{equation}
  \left.
  \begin{array}{rrcl} 
    \min & c(X) \\
    \forall b\in B &  \quad \underline{\mathscr{S}}_{b}\
      \leq\ \sum_{a\in B} X_{ba} \conj{\mathcal{Y}_{ba}}
      &\leq& \overline{\mathscr{S}}_{b}  \\ [0.4em]
    \forall b\in B & \underline{V}_b^2\ \le\ X_{bb} 
      &\le& \overline{V}_b^2 \\ [0.4em]
    \forall (b,a,1)\in L &
        X_{bb}X_{aa}  &\geq& |X_{ba}|^2 \\ [0.4em]
        & \multicolumn{3}{r}{X\in\mathbb{C}^{n\times n}\text{ hermitian}},
  \end{array}
  \right\} \label{socp-acopf}
\end{equation}
where the variable $X$ represent the terms found in the matrix $V\hermitian{V}$ (see first equality in Eq.~\eqref{lavpf}, also see Sect.~\ref{s:sdprel}), $c(X)$ is one the typical cost function and the constraint $X_{bb}X_{aa} \geq |X_{ba}|^2$ replaces the constraint that requires that the submatrix
${\small\left[\begin{array}{cc} X_{bb} &  X_{ba}\\ X_{ab} &  X_{aa}\end{array}\right]}$ is PSD (see Remark 1 in \cite{bose2015} in regards to dealing with additional constraints, such as, bounds on power flows). The numerical results presented in \cite{zuluaga-acopf} show the equivalence between the different problems stated above as well as trade-off between the time reduction resulting from the use of DD approximations and their precision.

\subsection{Improving the optimality gap of SDP relaxations}


\LEO{Some interesting improvements to optimality gap given by SDP relaxations and the corresponding exact formulations are proposed} in \cite{hijazi-globalopt}. \LEO{The SDP relaxation is approximated by a} sequence of formulations, called the ``determinant hierarchy'', introduced in \cite{hijazi-sdpcuts}. \LEO{These approximations are applied to each level of Lasserre's hierarchy.} Instead of requiring $U \succeq_{\mathsf{H}} 0$ for a given matrix in the SDP relaxation, the level $k$ of the determinant hierarchy imposes that all square sub-matrices of $U$ of size $\le k$ should have non-negative determinant, which yields polynomial constraints of degree $k$.

It is shown in \cite{madani2016,hijazi-sdpcuts} that the higher level of the determinant hierarchy (the level equal to the size of the matrix dimension) applied to the complex SDP relaxation Eq.~\eqref{Xsdp} corresponds to the replacement of constraint $X \succeq_{\mathsf{H}} 0$ by the constraints
\begin{equation}
  \forall {\cal S}\subseteq {\cal C} \in {\cal T}(G)\qquad \det(X_{\cal S})\ \ge\ 0, \label{treedec}
\end{equation}
where ${\cal T}(G)$ is a tree-decomposition \cite{diestel-minors} of the network $G$, ${\cal C}$ is one of the nodes of this tree-decomposition (which corresponds to a clique of $G$, thus, ${\cal C}\subseteq B$), and ${\cal S}$ is a sub-clique of ${\cal C}$. The matrix $X_{\cal S}$ is the submatrix of $X$ resulting from restricting $X$ to the columns and rows of $X$ corresponding to the buses in ${\cal S}$. Since real-life networks are usually sparse, it is possible to describe them with a tree-decomposition where the cliques have small size and, therefore, the number of sub-cliques is polynomial on the size of $B$.

\LEO{SDP relaxations and Reformulation-Linearization Technique (RLT) \cite{sh_alam} cuts were shown to be complementary in \cite{sdprltjogo}. This idea is applied to current and power expressions in \cite{hijazi-globalopt}.} Specifically, if Eqs.~\eqref{powercurrent}, \eqref{ohmlaw1} and \eqref{ohmlaw2a} are multiplied by their conjugate, we obtain the constraints
\begin{eqnarray*}
  \forall (b,a,h)\in L \qquad S_{bah}\conj{S_{bah}} &=& V_b\conj{V_b}I_{bah}\conj{I_{bah}} \\
  \forall (b,a,h)\in L_0\qquad I_{bah}\conj{I_{bah}} &=& |Y^{\mathsf{ff}}_{bah}|^2 V_b\conj{V_b} + Y^{\mathsf{ff}}_{bah} \conj{Y^{\mathsf{ft}}_{bah}}V_b\conj{V_a}+\\ && +\, Y^{\mathsf{ft}}_{bah}\conj{Y^{\mathsf{ff}}_{bah}} V_a\conj{V_b} + |Y^{\mathsf{ft}}_{bah}|^2 V_a\conj{V_a}\\
  \forall (b,a,h)\in L_1 \qquad I_{bah}\conj{I_{bah}} &=& |Y_{abh}^{\mathsf{tf}}|^2 V_a\conj{V_a} + Y_{abh}^{\mathsf{tf}}\conj{Y_{abh}^{\mathsf{tt}}} V_a\conj{V_b}+\\ && +\, Y_{abh}^{\mathsf{tt}}\conj{Y_{abh}^{\mathsf{tf}}} V_b\conj{V_a} + |Y_{abh}^{\mathsf{tt}}|^2 V_b\conj{V_b},
\end{eqnarray*}
which are then lifted to linear constraints using $X$ and the new variables $\widehat{\re{S_{bah}}}$, $\widehat{\im{S_{bah}}}$,  $\widehat{I_{bah}}$ and $\widehat{W\!I_{bah}}$ to represent the terms $(\re{S_{bah}})^2$, $(\im{S_{bah}})^2$, $|I_{bah}|^2$ and $W_{bb}\widehat{I_{bah}}$, respectively:
\begin{eqnarray}
  \forall (b,a,h)\in L \qquad \widehat{\re{S_{bah}}} &+& \widehat{\im{S_{bah}}}  = \widehat{W\!I_{bah}} \label{lift1}\\
  \forall (b,a,h)\in L_0\qquad \widehat{I_{bah}} &=& |Y^{\mathsf{ff}}_{bah}|^2 X_{bb} + Y^{\mathsf{ff}}_{bah} \conj{Y^{\mathsf{ft}}_{bah}}X_{ba}+ \label{lift2}\\ && +\, Y^{\mathsf{ft}}_{bah}\conj{Y^{\mathsf{ff}}_{bah}} X_{ab} + |Y^{\mathsf{ft}}_{bah}|^2 X_{aa} \notag\\
  \forall (b,a,h)\in L_1 \qquad \widehat{I_{bah}} &=& |Y_{abh}^{\mathsf{tf}}|^2 X_{aa} + Y_{abh}^{\mathsf{tf}}\conj{Y_{abh}^{\mathsf{tt}}} X_{ab}+ \label{lift3}\\ && +\, Y_{abh}^{\mathsf{tt}}\conj{Y_{abh}^{\mathsf{tf}}} X_{ba} + |Y_{abh}^{\mathsf{tt}}|^2 X_{bb}. \notag
\end{eqnarray}
More precisely, a different convex relaxation is obtained from Eq.~\eqref{Xsdp} with the following changes:
\begin{enumerate}
\item Eq.~\eqref{treedec} instead of constraint $X\succeq_{\mathsf{H}} 0$,
\item $\widehat{\re{S_{bah}}} + \widehat{\im{S_{bah}}}\le \overline{S}_{bah}^2$ as apparent flow limit at every arc $(b,a,h)\in L$,
\item inclusion of Eq.~\eqref{lift1}-\eqref{lift3} for each arc $(b,a,h)\in L$,
\item $0\le \widehat{I_{bah}}\le \overline{I}_{bah}^2$ as bounds for the squared current magnitude at each arc $(b,a,h)\in L$ (if the bound is not given, one could define $\overline{I}_{bah}:=\overline{S}_{bah}/\underline{V}_b$, see Eq.~\eqref{Vcurrentbound1}-\eqref{Vcurrentbound2rel}),
\item McCormick convex envelopes \cite{mccormick} for the lifted variable $\widehat{W\!I_{bah}}$ of the bilinear term $W_{bb}\widehat{I_{bah}}$:
\begin{align*}
  \forall(b,a,h)\in L\ \quad\ \underline{V}_b^2\widehat{I_{bah}}\ &\le\ \widehat{W\!I_{bah}}\ \le\ \overline{V}_b^2\widehat{I_{bah}}\\
  \forall(b,a,h)\in L\qquad \widehat{W\!I_{bah}}\ &\ge\ \overline{V}_b^2\widehat{I_{bah}} + W_{bb}\overline{I}_{bah}^2 - \overline{V}_b^2\overline{I}_{bah}^2\\
  \forall(b,a,h)\in L\qquad \widehat{W\!I_{bah}}\ &\le\ W_{bb}\overline{I}_{bah}^2 + \underline{V}_b^2\widehat{I_{bah}} - \underline{V}_{b}^2\overline{I}_{bah}^2,
\end{align*}
\item convexification of equations $\widehat{\re{S_{bah}}}=(\re{S_{bah}})^2$ and $\widehat{\im{S_{bah}}}=(\im{S_{bah}})^2$ defined by the inequalities:
\begin{align*}
  \forall(b,a,h)\in L\qquad (\re{S_{bah}})^2\ \le\ &\widehat{\re{S_{bah}}}\ \le\ \re{S_{bah}}(\underline{P}_{bah}+\overline{P}_{bah})-\underline{P}_{bah}\overline{P}_{bah}\\
  \forall(b,a,h)\in L\qquad (\im{S_{bah}})^2\ \le\ &\widehat{\im{S_{bah}}}\ \le\ \im{S_{bah}}(\underline{Q}_{bah}+\overline{Q}_{bah})-\underline{Q}_{bah}\overline{Q}_{bah},
\end{align*}
where $\underline{P}_{bah},\overline{P}_{bah},\underline{Q}_{bah},\overline{Q}_{bah}$ are lower/upper bounds for active and reactive power flows, one could take $-\underline{P}_{bah}=-\underline{Q}_{bah}=\overline{P}_{bah}=\overline{Q}_{bah}=\overline{S}_{bah}$.
\end{enumerate}
It turns out that this convex relaxation may be less tight w.r.t.~Eq.~\eqref{Xsdp}, but it can also be solved more efficiently, and hence larger instances can be tackled. See \cite{hijazi-envelopes} for more details.

An algorithm called SDP-based Bound Tightening based on Optimimality-based Bound Tightening \cite{branchcontract,gleixner-obbt} is proposed in \cite{hijazi-globalopt}. This algorithm alternately solves two types of sub-problems: (i) the convex ACOPF relaxation described above and (ii) a problem that minimizes (maximizes) one of the variables of this convex relaxation, subject to the original cost function bounded by the cost of a feasible solution, to find a tighter lower (upper) bound for the variable. The algorithm iterates between problem (i) and a set of problems of type (ii) obtained by considering different variables and whether it minimizes or maximizes the selected variable, the authors remark this second stage of the iterations can be parallelized through the consideration of the different sub-problems of type (ii).

\LEO{The approach proposed in \cite{hijazi-globalopt} is supported by numerical experiments on cases from benchmark libraries \cite[v19.05]{pglib} and \cite[v0.3]{nesta}. These experiments compare the gap obtained at the end of the algorithm with the gap of the two first levels of Lasserre's hierarchy. In the first library, the optimality gap obtained from the algorithm was less that 1\% on all tested instances (networks with 300 buses or less).}

\section{Conclusion}
\LEO{In this technical survey, we reviewed the alternating current optimal power flow problem and its modelling by mathematical programming. We presented continuous variable formulations in complex and real numbers, involving both polynomial and trigonometric terms, as well as relaxations based on several techniques, e.g.~semidefinite and second-order cone programming. Most of the formulations discussed in this survey have been modelled and tested on a few instances to verify consistency and reduce the occurrences of typos and errors (see the appendices below).

There are currently multiple challenges in this field. The first is to bridge the technical language of the power engineering community (both industrial and academic) to other complementary fields of knowledge: in this sense, this survey is intended as a contribution in this sense: we hope this survey will serve as a technical key to help operations researchers understand the details of power flow formulations in alternating current. The second challenge is to solve instances of this problem at national levels, i.e.~of very large size, to global optimality in relatively short times; there is a widespread belief (for obvious reasons) that the answer will come from relaxations, which is why we reviewed the latest contributions in Sect.~\ref{s:rel} --- some researchers, however, also think that good, special-purpose, fast solvers deployed on the original problem are essential to this purpose \cite{josz}. The third challenge, which we barely touched on, is to use the formulations presented above as a basis for more complicated formulations addressing security, design, distribution integration, and whole supply chain issues: we believe that such derived applications should be of interest to the operations research community at large.}

\ifspringer
\begin{acknowledgements}
\else
\section*{Acknowledgements}
\fi
\LEO{We are grateful to Cedric Josz for interesting technical discussions.}
\ifspringer
\end{acknowledgements}
\fi

{\small\noindent\textbf{Conflict of interest}\hspace*{4pt} The authors declare that they have no conflict of interest.}

\bibliographystyle{plain}
\bibliography{dr2,escobar}   

\begin{thebibliography}{10}

\bibitem{ahmadi2014}
A.~{Ahmadi} and A.~{Majumdar}.
\newblock {DSOS and SDSOS optimization: LP and SOCP-based alternatives to sum
  of squares optimization}.
\newblock In {\em 2014 48th Annual Conference on Information Sciences and
  Systems (CISS)}, pages 1--5, March 2014.

\bibitem{ahmadimajumdar}
A.~Ahmadi and A.~Majumdar.
\newblock {DSOS} and {SDSOS} optimization: {M}ore tractable alternatives to sum
  of squares and semidefinite optimization.
\newblock {\em SIAM Journal on Applied Algebra and Geometry}, 3(2):193--230,
  2019.

\bibitem{andersson}
G.~Andersson.
\newblock {Modelling and Analysis of Electric Power Systems}.
\newblock {\em EEH-Power Systems Laboratory, Swiss Federal Institute of
  Technology (ETH), Zürich, Switzerland}, 2008.

\bibitem{sdprltjogo}
K.~Anstreicher.
\newblock Semidefinite programming versus the reformulation-linearization
  technique for nonconvex quadratically constrained quadratic programming.
\newblock {\em Journal of Global Optimization}, 43:471--484, 2009.

\bibitem{pglib}
S.~{Babaeinejadsarookolaee}, A.~{Birchfield}, R.D. {Christie}, C.~{Coffrin},
  C.L. {DeMarco}, R.~{Diao}, M.~{Ferris}, S.~{Fliscounakis}, S.~{Greene},
  R.~{Huang}, C.~{Josz}, R.~{Korab}, B.C. {Lesieutre}, J.~{Maeght}, D.K.
  {Molzahn}, T.J. {Overbye}, P.~{Panciatici}, B.~{Park}, J.~{Snodgrass}, and
  R.~{Zimmerman}.
\newblock {The Power Grid Library for Benchmarking AC Optimal Power Flow
  Algorithms}.
\newblock Technical Report 1908.02788, arXiv, 2019.

\bibitem{beckopf}
A.~Beck, Y.~Beck, Y.~Levron, A.~Shtof, and L/~Tetruashvili.
\newblock Globally solving a class of optimal power flow problems in radial
  networks by tree reduction.
\newblock {\em Journal of Global Optimization}, 72:373--402, 2018.

\bibitem{fbbt-cocoa10}
P.~Belotti, S.~Cafieri, J.~Lee, and L.~Liberti.
\newblock Feasibility-based bounds tightening via fixed points.
\newblock In D.-Z. Du, P.~Pardalos, and B.~Thuraisingham, editors, {\em
  Combinatorial Optimization, Constraints and Applications (COCOA10)}, volume
  6508 of {\em LNCS}, pages 65--76, New York, 2010. Springer.

\bibitem{couenne}
P.~Belotti, J.~Lee, L.~Liberti, F.~Margot, and A.~W\"achter.
\newblock Branching and bounds tightening techniques for non-convex {MINLP}.
\newblock {\em Optimization Methods and Software}, 24(4):597--634, 2009.

\bibitem{bergen-book}
Arthur~R. Bergen and Vijay Vittal.
\newblock {\em {Power Systems Analysis}}.
\newblock Prentice Hall, Upper Saddle River, NJ, {Second} edition, 2000.

\bibitem{bienstock-acopf}
D.~Bienstock.
\newblock {\em Electrical {T}ransmission {S}ystem {C}ascades and
  {V}ulnerability: an {O}perations {R}esearch {V}iewpoint}.
\newblock Number~22 in MOS-SIAM Optimization. SIAM, Philadelphia, 2016.

\bibitem{bienstock2019}
D.~{Bienstock} and M.~{Escobar}.
\newblock {Stochastic Defense Against Complex Grid Attacks}.
\newblock {\em IEEE Transactions on Control of Network Systems}, 7(2):842--854,
  2020.

\bibitem{acopf_nphard_orl}
D.~Bienstock and A.~Verma.
\newblock Strong {NP}-hardness of {AC} power flows feasibility.
\newblock {\em Operations Research Letters}, 47:494--501, 2019.

\bibitem{bose2015}
S.~{Bose}, S.H. {Low}, T.~{Teeraratkul}, and B.~{Hassibi}.
\newblock {Equivalent Relaxations of Optimal Power Flow}.
\newblock {\em IEEE Transactions on Automatic Control}, 60(3):729--742, March
  2015.

\bibitem{oneill1}
M.~Cain, R.~O'Neill, and A.~Castillo.
\newblock History of optimal power flow and formulations.
\newblock Technical Report Staff Paper, Federal Energy Regulatory Commission,
  December 2012.

\bibitem{capitanescu2007}
F.~Capitanescu, M.~Glavic, D.~Ernst, and L.~Wehenkel.
\newblock Interior-point based algorithms for the solution of optimal power
  flow problems.
\newblock {\em Electric Power Systems Research}, 77(5):508 -- 517, April 2007.

\bibitem{Carpentier1962}
J.~Carpentier.
\newblock {Contribution \'a l'\'etude du dispatching \'economique}.
\newblock {\em Bulletin de la Soci\'et\'e Fran\c{c}aise des \'Electriciens},
  8(3):431--447, 1962.

\bibitem{carpentier1979}
J.~Carpentier.
\newblock Optimal power flows.
\newblock {\em International Journal of Electrical Power \& Energy Systems},
  1(1):3 -- 15, 1979.

\bibitem{chaojun2015}
G.~{Chaojun}, P.~{Jirutitijaroen}, and M.~{Motani}.
\newblock {Detecting False Data Injection Attacks in AC State Estimation}.
\newblock {\em IEEE Transactions on Smart Grid}, 6(5):2476--2483, September
  2015.

\bibitem{nesta}
C.~{Coffrin}, D.~{Gordon}, and P.~{Scott}.
\newblock {NESTA, The NICTA Energy System Test Case Archive}.
\newblock Technical Report 1411.0359, arXiv, 2014.

\bibitem{coffrin2016}
C.~{Coffrin}, H.~{Hijazi}, and P.~{Van Hentenryck}.
\newblock {The QC Relaxation: A Theoretical and Computational Study on Optimal
  Power Flow}.
\newblock {\em IEEE Transactions on Power Systems}, 31(4):3008--3018, July
  2016.

\bibitem{hijazi-envelopes}
C.~{Coffrin}, H.~{Hijazi}, and P.~{Van Hentenryck}.
\newblock {Strengthening the SDP Relaxation of AC Power Flows With Convex
  Envelopes, Bound Tightening, and Valid Inequalities}.
\newblock {\em IEEE Transactions on Power Systems}, 32(5):3549--3558, September
  2017.

\bibitem{ipopt}
COIN-OR.
\newblock {\em Introduction to IPOPT: A tutorial for downloading, installing,
  and using IPOPT}, 2006.

\bibitem{cvxpy}
S.~Diamond and S.~Boyd.
\newblock {CVXPY}: {A} {P}ython-embedded modeling language for convex
  optimization.
\newblock {\em Journal of Machine Learning Research}, 17:1--5, 2016.

\bibitem{diestel-minors}
R.~{Diestel}.
\newblock {\em Graph Minors}, pages 347--391.
\newblock Springer Berlin Heidelberg, Berlin, Heidelberg, 2017.

\bibitem{ecos}
A.~Domahidi, E.~Chu, and S.~Boyd.
\newblock {ECOS}: {A}n {SOCP} solver for embedded systems.
\newblock In {\em Proceedings of European Control Conference}, ECC, Piscataway,
  2013. IEEE.

\bibitem{ampl}
R.~Fourer and D.~Gay.
\newblock {\em The {AMPL} Book}.
\newblock Duxbury Press, Pacific Grove, 2002.

\bibitem{rebennack1}
S.~Frank, I.~Steponavice, and S.~Rebennack.
\newblock Optimal power flow: {A} bibliographic survey {I}. {F}ormulations and
  deterministic methods.
\newblock {\em Energy Systems}, 3:221--258, 2012.

\bibitem{rebennack2}
S.~Frank, I.~Steponavice, and S.~Rebennack.
\newblock Optimal power flow: {A} bibliographic survey {II}.
  {N}on-deterministic and hybrid methods.
\newblock {\em Energy Systems}, 3:259--289, 2012.

\bibitem{ghaddar2016}
B.~{Ghaddar}, J.~{Marecek}, and M.~{Mevissen}.
\newblock {Optimal Power Flow as a Polynomial Optimization Problem}.
\newblock {\em IEEE Transactions on Power Systems}, 31(1):539--546, January
  2016.

\bibitem{josz}
J.-C. Gilbert and C.~Josz.
\newblock Plea for a semidefinite optimization solver in complex numbers.
\newblock Technical Report hal-01422932, HAL Archives-Ouvertes, 2017.

\bibitem{snopt}
P.E. Gill.
\newblock {\em User's Guide for SNOPT 5.3}.
\newblock Systems Optimization Laboratory, Department of EESOR, Stanford
  University, California, February 1999.

\bibitem{gleixner-obbt}
A.M. {Gleixner}, T.~{Berthold}, B.~{M\"uller}, and S.~{Weltge}.
\newblock {Three enhancements for optimization-based bound tightening}.
\newblock {\em Journal of Global Optimization}, 67(4):731--757, April 2017.

\bibitem{overbye-book}
J.D. Glover, M.S. Sarma, and T.J. Overbye.
\newblock {\em {Power Systems Analysis and Design}}.
\newblock Cengage Learning, Stamford, CT, {Fourth} edition, 2008.

\bibitem{gomez1999}
A.~{G\'omez Exp\'osito} and E.~{Romero Ramos}.
\newblock Reliable load flow technique for radial distribution networks.
\newblock {\em IEEE Transactions on Power Systems}, 14(3):1063--1069, August
  1999.

\bibitem{hijazi-globalopt}
S.~{Gopinath}, H.~{Hijazi}, T.~{Wei{\ss}er}, H.~{Nagarajan}, M.~{Yetkin},
  K.~{Sundar}, and R.~{Bent}.
\newblock {Proving Global Optimality of ACOPF Solutions}.
\newblock In {\em 2020 Power Systems Computation Conference}, pages 1--6, June
  2020.

\bibitem{hijazi-sdpcuts}
H.~{Hijazi}, C.~{Coffrin}, and P.~{Van Hentenryck}.
\newblock {Polynomial SDP Cuts for Optimal Power Flow}.
\newblock In {\em 2016 Power Systems Computation Conference}, pages 1--7, June
  2016.

\bibitem{wei1998}
{Hua Wei}, H.~{Sasaki}, J.~{Kubokawa}, and R.~{Yokoyama}.
\newblock {An Interior Point Nonlinear Programming for Optimal Power Flow
  Problems with A Novel Data Structure}.
\newblock {\em IEEE Transactions on Power Systems}, 13(3):870--877, August
  1998.

\bibitem{huneault}
M.~Huneault and F.~Galiana.
\newblock A survey of the optimal power flow literature.
\newblock {\em IEEE Transactions on Power Systems}, 6(2):762--770, 1991.

\bibitem{jabr}
R.~Jabr.
\newblock Radial distribution load flow using conic programming.
\newblock {\em {IEEE} Transactions on Power Systems}, 21(3):1458--1459, 2006.

\bibitem{jabr2007}
R.~{Jabr}.
\newblock {A Conic Quadratic Format for the Load Flow Equations of Meshed
  Networks}.
\newblock {\em IEEE Transactions on Power Systems}, 22(4):2285--2286, November
  2007.

\bibitem{jabr2008}
R.~{Jabr}.
\newblock {Optimal Power Flow Using an Extended Conic Quadratic Formulation}.
\newblock {\em IEEE Transactions on Power Systems}, 23(3):1000--1008, August
  2008.

\bibitem{jabr2012}
R.~{Jabr}.
\newblock {Exploiting Sparsity in SDP Relaxations of the OPF Problem}.
\newblock {\em IEEE Transactions on Power Systems}, 27(2):1138--1139, May 2012.

\bibitem{kocuk2018}
B.~Kocuk, S.~Dey, and X.A. Sun.
\newblock {Matrix minor reformulation and SOCP-based spatial branch-and-cut
  method for the AC optimal power flow problem}.
\newblock {\em Mathematical Programming Computation}, 10(4):557--569, December
  2018.

\bibitem{zuluaga-acopf}
X.~Kuang, B.~Ghaddar, J.~Naoum-Sawaya, and L.~Zuluaga.
\newblock Alternative {LP} and {SOCP} hierarchies for {ACOPF} problems.
\newblock {\em IEEE Transactions on Power Systems}, 32(4):2828--2836, 2016.

\bibitem{lasserre}
J.~Lasserre.
\newblock Moments and sums of squares for polynomial optimization and related
  problems.
\newblock {\em Journal of Global Optimization}, 45:39--61, 2009.

\bibitem{lavaei}
J.~Lavaei and S.~Low.
\newblock Zero duality gap in optimal power flow problem.
\newblock {\em IEEE Transactions in Power Systems}, 27(1):92--107, 2012.

\bibitem{lesieutre2011}
B.C. {Lesieutre}, D.K. {Molzahn}, A.R. {Borden}, and C.L. {DeMarco}.
\newblock Examining the limits of the application of semidefinite programming
  to power flow problems.
\newblock In {\em 2011 49th Annual Allerton Conference on Communication,
  Control, and Computing (Allerton)}, pages 1492--1499, September 2011.

\bibitem{liu2017}
X.~{Liu} and Z.~{Li}.
\newblock {False Data Attacks Against AC State Estimation With Incomplete
  Network Information}.
\newblock {\em IEEE Transactions on Smart Grid}, 8(5):2239--2248, September
  2017.

\bibitem{yalmip}
J.~L\"ofberg.
\newblock {YALMIP}: {A} toolbox for modeling and optimization in {MATLAB}.
\newblock In {\em Proceedings of the International Symposium of Computer-Aided
  Control Systems Design}, volume~1 of {\em CACSD}, Piscataway, 2004. IEEE.

\bibitem{low2013}
S.H. {Low}.
\newblock {Convex Relaxation of Optimal Power Flow: A tutorial}.
\newblock In {\em 2013 IREP Symposium Bulk Power System Dynamics and Control -
  IX Optimization, Security and Control of the Emerging Power Grid}, pages
  1--15, August 2013.

\bibitem{low2014-1}
S.H. {Low}.
\newblock {Convex Relaxation of Optimal Power Flow—Part I: Formulations and
  Equivalence}.
\newblock {\em IEEE Transactions on Control of Network Systems}, 1(1):15--27,
  March 2014.

\bibitem{low2014-2}
S.H. {Low}.
\newblock {Convex Relaxation of Optimal Power Flow—Part II: Exactness}.
\newblock {\em IEEE Transactions on Control of Network Systems}, 1(2):177--189,
  June 2014.

\bibitem{madani2016}
R.~{Madani}, M.~{Ashraphijuo}, and J.~{Lavaei}.
\newblock {Promises of Conic Relaxation for Contingency-Constrained Optimal
  Power Flow Problem}.
\newblock {\em IEEE Transactions on Power Systems}, 31(2):1297--1307, March
  2016.

\bibitem{matlab2017a}
The MathWorks, Inc., Natick, MA.
\newblock {\em MATLAB R2017a}, 2017.

\bibitem{mccormick}
G.P. McCormick.
\newblock Computability of global solutions to factorable nonconvex programs:
  Part {I} --- {C}onvex underestimating problems.
\newblock {\em Mathematical Programming}, 10:146--175, 1976.

\bibitem{messine}
F.~Messine.
\newblock {\em M\'ethodes d'optimisation globale bas\'ees sur l'analyse
  d'intervalle pour la r\'esolution de probl\`emes avec contraintes (in
  {F}rench)}.
\newblock PhD thesis, Institut National Polytechnique de Toulouse, 1997.

\bibitem{sympy}
A.~Meurer, C.~Smith, M.~Paprocki, O.~\v{C}ert\'{i}k, S.~Kirpichev, M.~Rocklin,
  A.~Kumar, S.~Ivanov, J.~Moore, S~Singh, T.~Rathnayake, S.~Vig, B.~Granger,
  R.~Muller, F.~Bonazzi, H.~Gupta, S.~Vats, F.~Johansson, F.~Pedregosa,
  M.~Curry, A.~Terrel, \v{S}. Rou\v{c}ka, A.~Saboo, I.~Fernando, S.~Kulal,
  R.~Cimrman, and A.~Scopatz.
\newblock {SymPy}: {S}ymbolic computing in {P}ython.
\newblock {\em PeerJ Computer Science}, 3:e103, 2017.

\bibitem{molzahn2014-1}
D.K. {Molzahn} and I.A. {Hiskens}.
\newblock {Moment-Based Relaxation of the Optimal Power Flow Problem}.
\newblock In {\em 2014 Power Systems Computation Conference}, pages 1--7,
  August 2014.

\bibitem{molzahn2014-2}
D.K. {Molzahn} and I.A. {Hiskens}.
\newblock {Sparsity-Exploiting Moment-Based Relaxations of the Optimal Power
  Flow Problem}.
\newblock {\em IEEE Transactions on Power Systems}, 30(6):3168--3180, November
  2015.

\bibitem{molzahnsurvey}
D.K. {Molzahn} and I.A. {Hiskens}.
\newblock {A Survey of Relaxations and Approximations of the Power Flow
  Equations}.
\newblock {\em Foundations and Trends® in Electric Energy Systems},
  4(1-2):1--221, 2019.

\bibitem{molzahn2013}
D.K. {Molzahn}, J.T. {Holzer}, B.C. {Lesieutre}, and C.L. {DeMarco}.
\newblock {Implementation of a Large-Scale Optimal Power Flow Solver Based on
  Semidefinite Programming}.
\newblock {\em IEEE Transactions on Power Systems}, 28(4):3987--3998, November
  2013.

\bibitem{molzahn2014-3}
D.K. {Molzahn}, B.C. {Lesieutre}, and C.L. {DeMarco}.
\newblock {A Sufficient Condition for Global Optimality of Solutions to the
  Optimal Power Flow Problem}.
\newblock {\em IEEE Transactions on Power Systems}, 29(2):978--979, March 2014.

\bibitem{survey1993-1}
J.A. {Momoh}, M.E. {El-Hawary}, and R.~{Adapa}.
\newblock {A Review of Selected Optimal Power Flow Literature to 1993. I.
  Nonlinear and Quadratic Programming Approaches}.
\newblock {\em IEEE Transactions on Power Systems}, 14(1):96--104, February
  1999.

\bibitem{survey1993-2}
J.A. {Momoh}, M.E. {El-Hawary}, and R.~{Adapa}.
\newblock {A Review of Selected Optimal Power Flow Literature to 1993. II.
  Newton, Linear Programming and Interior Point Methods}.
\newblock {\em IEEE Transactions on Power Systems}, 14(1):105--111, February
  1999.

\bibitem{Monticelli1999}
A.~Monticelli.
\newblock {\em State Estimation in Electric Power Systems}.
\newblock Springer, Boston, 1999.

\bibitem{mosek8}
Mosek ApS.
\newblock {\em The \texttt{mosek} manual, Version 8}, 2016.

\bibitem{putinar2011}
M.~Putinar.
\newblock {Jean {B}ernard {L}asserre: {M}oments, Positive Polynomials and Their
  Applications (book review)}.
\newblock {\em Foundations of Computational Mathematics}, 11(4):489--497,
  August 2011.

\bibitem{rider2004}
M.J. {Rider}, V.L. {Paucar}, and A.V. {Garcia}.
\newblock Enhanced higher-order interior-point method to minimise active power
  losses in electric energy systems.
\newblock {\em IEE Proceedings - Generation, Transmission and Distribution},
  151(4):517--525, July 2004.

\bibitem{panciatici-acopf}
M.~Ruiz, J.~Maeght, A.~Mari\'e, P.~Panciatici, and A.~Renaud.
\newblock A progressive method to solve large-scale {AC} optimal power flow
  with discrete variables and control of the feasibility.
\newblock In {\em Proceedings of the Power Systems Computation Conference},
  volume~18 of {\em PSCC}, Piscataway, 2014. IEEE.

\bibitem{baron}
N.V. Sahinidis and M.~Tawarmalani.
\newblock {\em BARON 7.2.5: Global Optimization of Mixed-Integer Nonlinear
  Programs, {\em User's Manual}}, 2005.

\bibitem{salgado4}
E.~Salgado, C.~Gentile, and L.~Liberti.
\newblock Perspective cuts for the acopf with generators.
\newblock In P.~Daniele and L.~Scrimali, editors, {\em New trends in emerging
  complex real-life problems}, volume~1 of {\em AIRO Springer Series}.
  Springer, New York, 2018.

\bibitem{salgado3}
E.~Salgado, A.~Scozzari, F.~Tardella, and L.~Liberti.
\newblock Alternating current optimal power flow with generator selection.
\newblock In J.~Lee, G.~Rinaldi, and R.~Mahjoub, editors, {\em Combinatorial
  Optimization (Proceedings of ISCO 2018)}, volume 10856 of {\em LNCS}, pages
  364--375, 2018.

\bibitem{shectman}
J.P. Shectman and N.V. Sahinidis.
\newblock A finite algorithm for global minimization of separable concave
  programs.
\newblock {\em Journal of Global Optimization}, 12:1--36, 1998.

\bibitem{sh_alam}
H.D. Sherali and A.~Alameddine.
\newblock A new reformulation-linearization technique for bilinear programming
  problems.
\newblock {\em Journal of Global Optimization}, 2:379--410, 1992.

\bibitem{soltan2017}
S.~{Soltan} and G.~{Zussman}.
\newblock {Power grid state estimation after a cyber-physical attack under the
  AC power flow model}.
\newblock In {\em 2017 IEEE Power Energy Society General Meeting}, pages 1--5,
  July 2017.

\bibitem{sridhar2012}
S.~{Sridhar}, A.~{Hahn}, and M.~{Govindarasu}.
\newblock {Cyber–Physical System Security for the Electric Power Grid}.
\newblock {\em Proceedings of the IEEE}, 100(1):210--224, January 2012.

\bibitem{tinney1967}
W.F. {Tinney} and C.E. {Hart}.
\newblock {Power Flow Solution by Newton's Method}.
\newblock {\em IEEE Transactions on Power Apparatus and Systems},
  PAS-86(11):1449--1460, November 1967.

\bibitem{torres1998}
G.L. {Torres} and V.H. {Quintana}.
\newblock {An interior-point method for nonlinear optimal power flow using
  voltage rectangular coordinates}.
\newblock {\em IEEE Transactions on Power Systems}, 13(4):1211--1218, November
  1998.

\bibitem{VanNess1961}
J.E. {Van Ness} and J.H. {Griffin}.
\newblock {Elimination Methods for Load-Flow Studies}.
\newblock {\em Transactions of the American Institute of Electrical Engineers.
  Part III: Power Apparatus and Systems}, 80(3):299--302, April 1961.

\bibitem{python3}
G.~van Rossum and {\it et al.}
\newblock {\em Python Language Reference, version 3}.
\newblock Python Software Foundation, 2019.

\bibitem{wangpower}
H.~Wang, C.~Murillo-S\'anchez, R.~Zimmermann, and R.~Thomas.
\newblock On computational issues of market-based optimal power flow.
\newblock {\em IEEE Transactions on Power Systems}, 22(3):1185--1193, 2007.

\bibitem{branchcontract}
J.~M. Zamora and I.~E. Grossmann.
\newblock A branch and contract algorithm for problems with concave univariate,
  bilinear and linear fractional terms.
\newblock {\em Journal of Global Optimization}, 14:217:249, 1999.

\bibitem{matpower7}
R.~Zimmermann and C.~Murillo-S\'anchez.
\newblock {\em Matpower 7.0b1 {U}ser's {M}anual}.
\newblock Power Systems Engineering Research Center, 2018.

\bibitem{matpower}
R.~Zimmermann, C.~Murillo-Sanchez, and R.~Thomas.
\newblock {MATPOWER}: {S}teady-state operations, planning, and analysis tools
  for power systems research and education.
\newblock {\em IEEE Transactions on Power Systems}, 26(1):12--19, 2010.

\end{thebibliography}

\begin{appendices}

\section{Goal of the computational verification}
Many of the formulations discussed in this survey have been tested computationally. These tests were not designed to establish whether a formulation can be solved faster, or whether a relaxation is tighter, than another. The sheer complication of these formulations can be a formidable hurdle to their successful deployment, and previous experience from all of the authors confirmed that it is extremely difficult to remove all of the bugs. We therefore employed computational tests as a validity and consistency check.

\section{Modelling platforms}
We used three modelling platforms for implementing our formulations.
\begin{enumerate}
\item AMPL \cite{ampl} is a commercial, command-line interpreter which offers an incomparably elegant language, yielding code which is very similar to the formulations as they are presented mathematically on the written page. Its expression terms, objectives, and constraints may be quantified by indices varying on a set. AMPL has two serious limitations: (a) the amount of post-processing which can be expressed by its imperative sublanguage is limited (e.g.~there is no function for computing eigenvalues/eigenvectors or the inverse of a matrix); (b) there is no interface with most SDP solvers.
\item Python \cite{python3} is a {\it de facto} standard in ``scripting programming''. It is an interpreted language, with relatively low interpretation overhead CPU costs, and with a considerably large set of external modules (both interpreted and compiled), which allow the user to rapidly code almost anything. An equally large corpus of online documentation makes it possible to solve and issues using a simple internet query. We used the {\tt cvxpy} \cite{cvxpy} MP modelling interface, which allows the coding and solution of SDPs and SOCPs using a range of solvers.
\item \textsc{Matlab} \cite{matlab2017a} is a well-known commercial ``general-purpose'' applied mathematical software package. It offers very good MP capabilities through a range of interfaces. We used YALMIP \cite{yalmip}, which, notably, also connects to SDP solvers.
\end{enumerate}
Specifically, we implemented and tested:
\begin{itemize}
\item the real cartesian $(S,I,V)$-formulation (Sect.~\ref{s:sivreal}), the real cartesian voltage-only QCQP formulation (Sect.~\ref{s:vcart}), the real polar NLP formulation (also in Sect.~\ref{s:vcart}), and Jabr's (real) relaxation (Sect.~\ref{s:jabr}) using AMPL;
\item the complex SDP relaxation Eq.~\eqref{vsdp} (Sect.~\ref{s:sdprel}) using {\tt cvxpy} on Python3;
\item the real matrix formulation (Sect.~\ref{s:matrix}) using YALMIP on \textsc{Matlab}.
\end{itemize}

We solved our formulations with a variety of solvers, some global and some local (used within a multi-start heuristic): Baron \cite{baron}, Couenne \cite{couenne}, ECOS \cite{ecos}, Mosek \cite{mosek8}, IPOpt \cite{ipopt}, Snopt \cite{snopt}. We remark that Baron cannot deal with trigonometric functions. We found that {\tt cvxpy} was able to pass complex number SDPs to ECOS correctly, but not (always) to Mosek. The global solvers for NLP (Baron, Couenne), could never certify global optima, even for the smallest instances, testifying to the practical hardness of the ACOPF.

On the other hand, Baron's ``upper bounding heuristic'', consisting in a multi-start on various local NLP solvers, yielded best solutions. Our benchmark was formed by small instances in \textsc{Matlab}'s \textsc{MatPower}'s \cite{matpower7} {\tt data} folder, and our comparison stone by the optima found by \textsc{MatPower}'s own local NLP solver --- a \textsc{Matlab} implementation of a standard interior point method algorithm --- from the local optima stored in the instances. 

We also implemented a symbolic computation code (provided by {\tt sympy} \cite{sympy}) in order to derive the the real expressions of formulations in real numbers from the corresponding complex expressions occurring in their complex counterparts.

\section{AMPL code of the $(S,I,V)$-formulation}
In this section we list the AMPL code of the real $(S,I,V)$-formulation (Sect.~\ref{s:sivreal}). We do not report all our AMPL code for brevity. We think nonetheless that the code below will be sufficient to give the general idea (all our other code is available upon request --- contact LL). We start with declarations of constants, sets and parameters.

{\footnotesize
\seplin
\begin{verbatim}
### constants
param Inf := 1e30;
param Eps := 1e-6;
param Pi := 4*atan(1);
param myZero := 1e-9;

### sets

# max number of parallel branches in instance
param maxParBranches integer, >0, default 1;
set PB := 1..maxParBranches;

# buses
set B;

## lines (directed: edges correspond to antisymmetric arcs)
# the given edges (b,a): the transformer is on b
set L0 within {B,B,PB};  
# set of all arcs
set L default L0 union {(a,b,i) in {B,B,PB} : (b,a,i) in L0}; 
# set of all antiparallel arcs
set L1 default L diff L0;

# set of generators at node
set G{b in B} default {};

### parameters

# bus type (2=generator, 3=reference)
param busType{B} integer;

# cost coefficients (minimization of power generation costs)
param Kcard integer, >= 0, default 2;
set K := 0..Kcard;
# initialization: only linear terms in P (quadratic in V)
param C{b in B, g in G[b], k in K} default if k == 1 then 1 else 0;

# real power demand at buses (there can be nodes with negative demands)
param SDR{B} default 0; 
# reactive power demand at buses (appears in MATPOWER documentation)
param SDC{B} default 0;

# power bounds at generators (real, imaginary)
param SLR{b in B, G[b]} default -Inf;
param SLC{b in B, G[b]} default -Inf;
param SUR{b in B, g in G[b]} >= SLR[b,g], default Inf;
param SUC{b in B, g in G[b]} >= SLC[b,g], default Inf;

# upper power magnitude bounds on links - symmetric
param SU{L} >= 0, default Inf; 

# voltage magnitude bounds at buses
param VL{B} default 0; # can't have negative moduli
param VU{b in B} >= VL[b], default Inf;

# shunt parameters at buses
param shR{B} default 0;  # MatPower's Gs
param shC{B} default 0;  # MatPower's Bs

# status of a branch
param status{L} default 1;

# Y matrix (Ohm's law in AC) data, only defined on given arcs in L0
param r{L0} default 0;
param x{L0} default 0;
param bb{L0} default 0;
param tau{L0} default 1;
param nu{L0} default 0;  # translated to radians by mpc2dat.py

# phase difference bounds (only across lines in L0)
#   translated to radians by mpc2dat.py ##191026: not using these
param pdLB{L0} default -Pi;
param pdUB{L0} default Pi;

# the 2x2 complex Y matrix appearing in Ohm's law for a line (b,a) in L0
#  (initialized from r,x,bb,tau,nu in a .run file)
param YffR{L0} default 0;
param YffC{L0} default 0;
param YftR{L0} default 0;
param YftC{L0} default 0;
param YtfR{L0} default 0;
param YtfC{L0} default 0;
param YttR{L0} default 0;
param YttC{L0} default 0;
\end{verbatim}
\seplin
}

\noindent Next, we introduce the declarations for decision variables.

{\footnotesize
\seplin
\begin{verbatim}
### decision variables

# voltage (real, imaginary)
var VR{b in B} <= VU[b], >= -VU[b];  # real
var VC{b in B} <= VU[b], >= -VU[b];  # imaginary

# power generation (real, imaginary)
var SgenR{b in B, g in G[b]} >= SLR[b,g], <= SUR[b,g];
var SgenC{b in B, g in G[b]} >= SLC[b,g], <= SUC[b,g];

# V2 = |V|^2
var V2{b in B} >= VL[b]^2, <= VU[b]^2;

# current (real, imaginary)
var IR{L};
var IC{L};

# power injected on line at bus (real, imaginary)
var SR{L};
var SC{L};
\end{verbatim}
\seplin
}

\noindent Lastly, we detail the definitions of objective function and constraints.

{\footnotesize
\seplin
\begin{verbatim}
### objective function

# generation cost (WARNING: to remove quadratic power terms set C[g,2]=0)
minimize gencost: 
  sum{b in B, g in G[b]} 
    (C[b,g,2]*SgenR[b,g]^2 + C[b,g,1]*SgenR[b,g] + C[b,g,0]);

### constraints

# power flow (real, imaginary)
subject to powerflowR{b in B}:
   SDR[b] + sum{(b,a,i) in L} SR[b,a,i] = 
  -shR[b]*V2[b] + sum{g in G[b]} SgenR[b,g];
subject to powerflowC{b in B}:
   SDC[b] + sum{(b,a,i) in L} SC[b,a,i] = 
   shC[b]*V2[b] + sum{g in G[b]} SgenC[b,g];

# definition of power in function of current and voltage (real, imaginary)
subject to powerinjR{(b,a,i) in L}: 
  SR[b,a,i] = VR[b]*IR[b,a,i] + VC[b]*IC[b,a,i];
subject to powerinjC{(b,a,i) in L}: 
  SC[b,a,i] = VC[b]*IR[b,a,i] - VR[b]*IC[b,a,i];

# Ohm's law ((b,a):real,imaginary; (a,b):real,imaginary)
subject to ohm1R{(b,a,i) in L0}:
  IR[b,a,i] = YffR[b,a,i]*VR[b]-YffC[b,a,i]*VC[b] 
            + YftR[b,a,i]*VR[a]-YftC[b,a,i]*VC[a];
subject to ohm1C{(b,a,i) in L0}:
  IC[b,a,i] = YffR[b,a,i]*VC[b]+YffC[b,a,i]*VR[b] 
            + YftR[b,a,i]*VC[a]+YftC[b,a,i]*VR[a];
subject to ohm2R{(b,a,i) in L0}:
  IR[a,b,i] = YtfR[b,a,i]*VR[b]-YtfC[b,a,i]*VC[b] 
            + YttR[b,a,i]*VR[a]-YttC[b,a,i]*VC[a];
subject to ohm2C{(b,a,i) in L0}:
  IC[a,b,i] = YtfR[b,a,i]*VC[b]+YtfC[b,a,i]*VR[b] 
            + YttR[b,a,i]*VC[a]+YttC[b,a,i]*VR[a];

# power bound on lines b->a defined on I
subject to powerbound{(b,a,i) in L : SU[b,a,i]>0 and SU[b,a,i]<Inf}:
  SR[b,a,i]^2 + SC[b,a,i]^2 <= SU[b,a,i]^2;

# definition of V2
subject to V2def{b in B}: V2[b] = VR[b]^2 + VC[b]^2;

# bounds on phase difference
subject to phasediff1{(b,a,i) in L0}:
  VC[b]*VR[a] - VR[b]*VC[a] <= tan(pdUB[b,a,i]) * (VR[b]*VR[a] + VC[b]*VC[a]);
subject to phasediff2{(b,a,i) in L0}:
  VC[b]*VR[a] - VR[b]*VC[a] >= tan(pdLB[b,a,i]) * (VR[b]*VR[a] + VC[b]*VC[a]);
subject to phasediff3{(b,a,1) in L0}: VR[b]*VR[a] + VC[b]*VC[a] >= 0;

# reference bus: there had better be just one reference -- check in .run
subject to reference1{b in B : busType[b] == 3}: VC[b] = 0;
subject to reference2{b in B : busType[b] == 3}: VR[b] >= 0;
\end{verbatim}
\seplin
}

\noindent The {\tt .dat} files contain essentially the same information as the corresponding {\tt .m} \textsc{Matlab} files distributed with \textsc{MatPower}. For example, the smallest instance file {\tt case5.dat} is as follows:

{\footnotesize
\seplin
\begin{verbatim}
param maxParBranches := 1 ;
param : B : busType SDR SDC VL VU Vm Va shR shC :=
  1   2 0.0 0.0 0.9 1.1 1.0 0.0 0.0 0.0
  2   1 3.0 0.9861 0.9 1.1 1.0 0.0 0.0 0.0
  3   2 3.0 0.9861 0.9 1.1 1.0 0.0 0.0 0.0
  4   3 4.0 1.3147 0.9 1.1 1.0 0.0 0.0 0.0
  5   2 0.0 0.0 0.9 1.1 1.0 0.0 0.0 0.0 ;
set G[1] := 1 2 ;
set G[3] := 1 ;
set G[4] := 1 ;
set G[5] := 1 ;
param : SLR SLC SUR SUC :=
  1 1   0.0 -0.3 0.4 0.3
  1 2   0.0 -1.275 1.7 1.275
  3 1   0.0 -3.9 5.2 3.9
  4 1   0.0 -1.5 2.0 1.5
  5 1   0.0 -4.5 6.0 4.5 ;
param : L0 : status SU r x bb tau nu pdLB pdUB :=
  1 2 1   1 4.0 0.00281 0.0281 0.00712 1.0 0.0 -1.57079632679 1.57079632679
  1 4 1   1 1e+30 0.00304 0.0304 0.00658 1.0 0.0 -1.57079632679 1.57079632679
  1 5 1   1 1e+30 0.00064 0.0064 0.03126 1.0 0.0 -1.57079632679 1.57079632679
  2 3 1   1 1e+30 0.00108 0.0108 0.01852 1.0 0.0 -1.57079632679 1.57079632679
  3 4 1   1 1e+30 0.00297 0.0297 0.00674 1.0 0.0 -1.57079632679 1.57079632679
  4 5 1   1 2.4 0.00297 0.0297 0.00674 1.0 0.0 -1.57079632679 1.57079632679 ;
param Kcard := 2 ;
param C :=
  1 1 1   1400.0
  1 2 1   1500.0
  3 1 1   3000.0
  4 1 1   4000.0
  5 1 1   1000.0 ;
\end{verbatim}
\seplin
}

\end{appendices}

\end{document}